\documentclass[11pt, a4paper, reqno, oneside]{article}

\usepackage{amssymb,amsmath}
\DeclareUnicodeCharacter{0335}{}
\usepackage{tikz-cd}
\usepackage{authblk} 
\usepackage{physics} 
\usepackage{tikz}
\usepackage{titlesec} 
\usetikzlibrary{babel}
\usepackage{adjustbox}
\usepackage{url}
\usepackage{xcolor}
\usepackage{yhmath}
\usepackage{eso-pic}
\usepackage[T1]{fontenc}
\usepackage[utf8]{inputenc}
\usepackage[spanish,english]{babel}
\usepackage{csquotes}
\usepackage{pst-node}
\usepackage[english]{babel}
\usepackage{latexsym}
\usepackage{euscript}
\usepackage{mathtools}
\usepackage[all]{xy}
\usepackage{mathrsfs}
\usepackage[percent]{overpic}
\usepackage{setspace}
\usepackage{graphicx}
\usepackage{pdfsync}
\usepackage{amsthm}
\usepackage{thmtools}
\usepackage{libertine}
\usepackage{yfonts}
\usetikzlibrary{calc}
\usetikzlibrary{decorations.pathmorphing}

\tikzset{curve/.style={settings={#1},to path={(\tikztostart)
    .. controls ($(\tikztostart)!\pv{pos}!(\tikztotarget)!\pv{height}!270:(\tikztotarget)$)
    and ($(\tikztostart)!1-\pv{pos}!(\tikztotarget)!\pv{height}!270:(\tikztotarget)$)
    .. (\tikztotarget)\tikztonodes}},
    settings/.code={\tikzset{quiver/.cd,#1}
        \def\pv##1{\pgfkeysvalueof{/tikz/quiver/##1}}},
    quiver/.cd,pos/.initial=0.35,height/.initial=0}

\tikzset{tail reversed/.code={\pgfsetarrowsstart{tikzcd to}}}
\tikzset{2tail/.code={\pgfsetarrowsstart{Implies[reversed]}}}
\tikzset{2tail reversed/.code={\pgfsetarrowsstart{Implies}}}
\tikzset{no body/.style={/tikz/dash pattern=on 0 off 1mm}}

\usepackage[colorlinks=true,linktocpage=true, pagebackref=false, citecolor=black,linkcolor=black]{hyperref}
\normalbaselines
\makeatletter

\def\ps@myfancy{\let\@mkboth\markboth
 \def\@evenhead{\vbox{\hsize\textwidth 
 \hbox to \textwidth{\sf\mdseries\thepage 
 \rule[-.6ex]{0mm}{2mm} \hfill\sf\large\leftmark}
 \vskip 1pt \hrule}}
 \def\@oddhead{\vbox{\hsize\textwidth 
 \hbox to \textwidth{{\sf\large\leftmark}
 \rule[-.6ex]{0mm}{2mm} \hfill\sf\mdseries{\thepage}}
 \vskip 1pt \hrule}}}

\def\ps@myfancyplain{
 \def\@evenhead{\vbox{\hsize\textwidth%
 \rule[-.6ex]{0mm}{2mm} \hfill }
 \vskip 1pt \hrule
 \vskip\headsep
 \vskip\textheight
 \vskip1pc
 \hbox to \textwidth{\sf\mdseries\thepage 
 \rule[-6ex]{0mm}{2mm} \hfill }}
 \def\@oddhead{\vbox{\hsize\textwidth 
 \vskip 1pt\hrule
 \vskip\headsep
 \vskip\textheight
 \vskip2pc
 \hbox to \textwidth{\hfill\rule[.4ex]{1pc}{2.5pt}
 \sf\mdseries\thepage}
}}}

\def\ps@myemptyfun{
 \def\@evenhead{\vbox{\hsize\textwidth
 \rule[-.6ex]{0mm}{2mm} \hfill }
 \vskip 1pt 
 \vskip\headsep
 \vskip\textheight
 \vskip1pc
 \hbox to \textwidth{\sf\mdseries\thepage 
 \rule[-0.6ex]{0mm}{2mm} 
 \hfill }}
 \def\@oddhead{\vbox{\hsize\textwidth 
 \vskip 1pt 
 \vskip\headsep
 \vskip\textheight
 \vskip2pc
}}}

\makeatother
\providecommand{\proofname}{Demostraci\'on.}
 {\par\noindent{\it Demostraci\'on. }\nopagebreak\normalsize}%

 {\par\noindent{\it #1. }\nopagebreak\normalsize}%
 {\hfill\linebreak[2]\hspace*{\fill}$\square$\\[-1pt]}
\makeatother

\def\sqbullet{\raise.2ex\hbox{\vrule width 3.5pt height 3.5pt}}

{}

{\em}

\newcounter{substep}
\def\thesubstep{\arabic{substep}}

{\em}

\newcounter{subsubstep}
\def\thesubsubstep{\arabic{subsubstep}}

{\em}

\numberwithin{figure}{section}

{}

\newtheoremstyle{mystyle}
  {}
  {}
  {\itshape}
  {}
  {\sf \bfseries}
  {}
{ }
  {\thmname{#1}\thmnumber{{\textcolor{blue}{\, \hspace{-1mm}#2.}}}\thmnote{ (#3)}}

\theoremstyle{mystyle}
\definecolor{royalblue(web)}{rgb}{0.25, 0.41, 0.88}
\hypersetup{
    colorlinks=true,
    linkcolor=royalblue(web),
    filecolor=magenta,      
    urlcolor=cyan,
    pdftitle={Overleaf Example},
    pdfpagemode=FullScreen,
    }

\titleformat{\section}
  {\normalfont\LARGE\bfseries}{\thesection.}{1em}{}
\titleformat{\subsection}
  {\normalfont\Large\bfseries}{\thesubsection.}{1em}{}
\titleformat{\subsubsection}
  {\normalfont\normalsize\itshape}{\thesubsubsection.}{1em}{}

\newtheorem{Teor}{Theorem}[section]

\newtheorem{Prop}[Teor]{Proposition}
\newtheorem{Coro}[Teor]{Corollary}

\newtheorem{Defi}[Teor]{Definition}

\newtheorem{Lema}[Teor]{Lemma}

\newtheorem{Ejems}[Teor]{Examples}

\newtheorem{Obse}[Teor]{Remark}

\newcommand{\R}{{\mathbb R}}
 \newcommand{\C}{{\mathbb C}}

\newcommand{\Cc}{\mathcal{C}} 

 \newcommand{\B}{{\cal B}}

\newcommand{\mor}{\operatorname{Mor}}

\newcommand{\mail}[1]{\small\href{mailto:#1}{#1}}

\newenvironment{Abstract}
{
\begin{center}
\textbf{Abstract}\\
\vspace{0.25cm}
\begin{minipage}{14.5cm}}
{\footnotesize
\end{minipage}
\end{center}}
\begin{document}


	\begin{center}
		{\huge {\bfseries Bessel Potential Spaces and Complex Interpolation: Continuous embeddings}\par}
		\vspace{1cm}
		
\begin{tabular}{l@{\hskip 2cm}l} 
	{\Large José C. Bellido}{\small\textsuperscript{1}} & {\Large Guillermo García-Sáez}{\small\textsuperscript{1}} \\
	\mail{josecarlos.bellido@uclm.es} & \mail{guillermo.garciasaez@uclm.es}
\end{tabular}
\vspace{5mm}

\textsc{\textsuperscript{1}ETSII, Departamento de Matem\'aticas\\ Universidad de Castilla-La Mancha} \\
		Campus Universitario s/n, 13071 Ciudad Real, Spain. \\ \vspace{5mm}
\end{center}
\begin{Abstract}
Bessel potential spaces, introduced in the 1960s, are derived through complex interpolation between Lebesgue and Sobolev spaces, making them intermediate spaces of fractional differentiability order. Bessel potential spaces have recently gained attention due to their identification with the space associated to the Riesz fractional gradient. This paper explores Bessel potential spaces as complex interpolation spaces, providing original proofs of fundamental properties based on abstract interpolation theory. Main results include a direct proof of norm equivalence, continuous embeddings, and the relationship with Gagliardo spaces. 
\end{Abstract}

\noindent {\bf Keywords:} Bessel potential spaces, complex interpolation method, Sobolev spaces of fractional order, continuous embeddings. 

\noindent {\bf AMS Subject Classification:} 46B70

\tableofcontents

\section{Introduction}

Bessel potential spaces were introduced in the 1960s, as documented in \cite{Aronszajn1961}, driven by the study of partial differential equations (PDE), particularly in the linear case. Remarkably,  Bessel potential spaces are the image of trace operators acting on Sobolev spaces \cite{LionsMagenes1972}. These spaces can be derived through complex interpolation between Lebesgue and Sobolev spaces \cite{Calderon1963,Calderon1961,Lions1960,LionsMagenes1972}, positioning them as intermediate spaces between those spaces, or as Sobolev spaces of fractional differentiability order. Another class of fractional Sobolev spaces, known as Gagliardo or Sobolev-Slobodeckij spaces, is obtained via the real method of interpolation rather than the complex method. In recent decades, Gagliardo spaces have garnered significant attention due to their direct connection to the fractional Laplacian, or more generally, the fractional $p$-Laplacian. In contrast, Bessel potential spaces have received limited attention until recently.

Both Gagliardo and Bessel potential spaces are genuine fractional Sobolev spaces. However, the term "fractional Sobolev spaces" has become more commonly associated with Gagliardo spaces due to their extensive use in PDE analysis over the past three or four decades. Nonetheless, there has been a renewed interest in Bessel potential spaces following the pioneering works \cite{ShiehSpector2015,ShiehSpector2018}, which established a connection with the Riesz fractional gradient. The Riesz fractional gradient of order $s$ for a smooth compactly supported function $u$ is defined as
\[\nabla^s u(x)=c_{n,s} \int_{\mathbb{R}^n} \frac{u(x)-u(y)}{|x-y|^{n+s}}\frac{x-y}{|x-y|}\,dy,\]
where $c_{n,s}$ is a normalization constant dependent on the dimension $n$ and the fractional index $s$.

In \cite{ShiehSpector2015,ShiehSpector2018}, it was demonstrated that Bessel potential spaces coincide, with equivalence of norms, with the closure of $C_c^\infty(\mathbb R^n)$ with respect to the norm
\[ \|u\|=\|u\|_{L^p(\mathbb{R}^n)}+\|\nabla^s u\|_{L^p(\mathbb{R}^n,\mathbb{R}^n)}.\]
The Riesz fractional gradient exhibits several notable properties. In \cite{Silhavy2020}, the uniqueness of the Riesz fractional gradient, up to a multiplicative constant, was established under natural conditions (invariance under translations and rotations, homogeneity under dilations, and certain continuity in the sense of distributions). Furthermore, for a Sobolev function $u \in W^{1,p}(\mathbb{R}^n)$, $D^s u$ converges strongly to $Du$ in $L^p(\mathbb R^n)$ as $s\to 1^-$ \cite{BellidoCuetoMoraCorral2021}. Additionally, a fractional divergence of order $s$ can be defined to ensure integration by parts holds (see \cite{BellidoCuetoMoraCorral2020} and references therein). Moreover, fractional laplacian in $\mathbb{R}^n$ can be obtained as the fractional divergence of the Riesz fractional gradient. Consequently, the Riesz fractional gradient $D^s u$ is a truly fractional differential object, and this, combined with the identification of Bessel potential spaces in \cite{ShiehSpector2015}, makes Bessel potential spaces an intriguing and valuable framework for the analysis of fractional PDEs.

Applications of Bessel potential spaces include the study of a new family of fractional PDEs \cite{ShiehSpector2015,ShiehSpector2018}, variational inequalities involving the Riesz fractional gradient \cite{Campos2024,CamposRodrigues2023II,campos2023Primero,LoRodriges2023}, fractional elasticity problems \cite{BellidoCuetoMoraCorral2020,BellidoCuetoMoraCorral2021}, and the development of a new concept of fractional perimeter \cite{BrueCalziComiStefani2022,ComiStefani2023,ComiStefani2019}.

A significant factor contributing to the growing interest in Bessel potential spaces is their robust structural properties from the perspective of Functional Analysis, making them highly applicable and useful in the study of fractional PDEs. Specifically, \cite{ShiehSpector2015} established continuous embeddings of Bessel potential spaces into Lebesgue spaces, and \cite{ShiehSpector2018} demonstrated the compactness of these embeddings. These results fundamentally rely on the identification of the spaces via the Riesz fractional gradient, as previously mentioned. However, the fact that Bessel potential spaces are interpolation spaces has not been fully exploited in those analyses.
\newpage
Our objective here is to delve deeper into the nature of Bessel potential spaces as complex interpolation spaces to directly derive their main properties. While it is frequently mentioned in the literature that Bessel potential spaces are obtained through complex interpolation of Lebesgue and Sobolev spaces, most references, such as \cite[Theorem 7.63]{Adams1975}, provide this information without a proof. In this paper, we clarify the relationship between Bessel potential spaces and interpolation theory, offering original proofs of all fundamental facts. Our approach is based solely on well-established principles of abstract interpolation theory, without relying on results from other general families of spaces, such as Triebel-Lizorkin spaces \cite{Triebel1995,Triebel1992}. Consequently, our proofs are specific to Bessel potential spaces and novel in the literature, despite some results having been previously demonstrated using different techniques.

This work is the first in a series of three papers examining Bessel potential spaces from the perspective of complex interpolation theory. This installment focuses on continuous embeddings into Lebesgue and Hölder spaces and the relationship between Bessel potential spaces and Gagliardo spaces. The second paper will explore the relationship between Bessel spaces and the Riesz fractional gradient, addressing some ambiguities in the literature regarding the definition of spaces involving the Riesz fractional gradient, and examining the asymptotics as $s \to 1$ and $s \to 0$ using complex interpolation techniques. The third paper will address the compactness of the aforementioned embeddings, deriving compactness from complex interpolation theory results and providing alternative proofs based on the fractional gradient.

The structure of this paper is as follows: Section 2 provides a brief introduction to interpolation theory, covering both the real and complex methods of interpolation, including the main results utilized in this paper and the relationship between the two methods. Section 3 presents the main results of the paper. Notably, Theorem \ref{NormBessel} directly establishes that the definition of Bessel potential spaces coincides with norm equivalence with $[L^p(\mathbb{R}^n), W^{1,p}(\mathbb{R}^n)]_s$, $s\in[0,1]$, the complex interpolation of order $s$ between Lebesgue and Sobolev spaces. Theorem \ref{FSET} compiles the continuous embeddings of Bessel spaces in various cases—subcritical, critical, and supercritical—highlighting that our proofs are novel in the literature. Theorem \ref{Hilbertcase} establishes the equality of Bessel and Gagliardo spaces in the Hilbert case ($p=2$), while Theorems \ref{contiguity} and \ref{nesting} provide contiguity and nesting properties among Bessel and Gagliardo spaces as a consequence of the relationship between real and complex interpolation. Finally, Section 3.4 addresses the open question of whether the Bessel space $H^{s,1}$ ($p=1$) is an interpolation space.

\section{Preliminaries on Interpolation Theory}

This section is a brief introduction to Interpolation Theory, including the concepts and results we will use in the following. For a complete development of the main ideas, we refer to \cite{BerghLofstrom1976,GarciaSaez2024,Lunardi2018,Triebel1995} \\

To motivate the framework of Interpolation Theory, we will briefly look at what happens with Lebesgue spaces $L^p(\Omega,\mu)$, for $1\leq p\leq \infty$, where $(\Omega,\mu)$ is a measure space. Typically, when there is no risk of confusion we omit the measure $\mu$ in the notation, as for the Lebesgue measure for instance. Given $f\in L^p(\Omega,\mu)$ and $t>0$, we can write 
$$f=f\chi_{|f|>t}+f\chi_{|f|\leq t}=:f_1+f_2,$$ 
where clearly $f_2$ is a function such that $\norm{f_2}_\infty\leq t$. Further, ç$$\norm{f_1}_1=\int_\Omega |f_1|\,d\mu=\int_{|f|>t}|f|\,d\mu\leq \int_{|f|>t}|f|\frac{|f|^{p-1}}{t^{p-1}}\,d\mu=t^{1-p}\norm{f}_p^p.$$ 
These two estimates imply that for every $f\in L^p(\Omega,\mu)$, we can find $f_1\in L^1(\Omega,\mu)$ and $f_2\in L^\infty(\Omega,\mu)$ such that $f=f_1+f_2\in L^1(\Omega,\mu)+L^\infty(\Omega,\mu)$ and 
$$\norm{f_1}_1+\norm{f_2}_\infty\leq t^{1-p}\norm{f}_p^p+t.$$ 
Choosing $t=(p-1)^{1/p}\norm{f}_p$, it yields that $$\norm{f_1}_1+\norm{f_2}_\infty\leq \left((p-1)^{(1-p)/p}+(p-1)^{1/p}\right)\norm{f}_p,$$ hence 
$$\norm{f}_{L^1+L^\infty}:=\inf\left\{\norm{g}_1+\norm{h}_\infty: f=g+h, g\in L^1,h\in L^\infty\right\}\leq \left((p-1)^{(1-p)/p}+(p-1)^{1/p}\right)\norm{f}_p.$$ This implies that for every $1<p<\infty$, $L^p(\Omega,\mu)\xhookrightarrow{}L^1(\Omega,\mu)+L^\infty(\Omega,\mu)$, where $\xhookrightarrow{}$ denotes the usual continuous embedding between normed spaces.
Now, we take $1\leq p_0<p_1\leq \infty$, and consider 
$$\frac{1}{p}:=\frac{1-\theta}{p_0}+\frac{\theta}{p_1},\,\theta\in (0,1).$$
By H\"older's inequality \begin{align*}
    \norm{f}_p^p=\int_\Omega |f|^p\,d\mu=\int_\Omega |f|^{(1-\theta)p}|f|^{\theta p}\,d\mu\leq \left(\int_\Omega |f|^{p_0}\,d\mu\right)^{p(1-\theta)/p_0}\left(\int_\Omega |f|^{p_1}\,d\mu\right)^{p\theta /p_1},
\end{align*}
for all $f\in L^{p_0}(\Omega,\mu)\cap L^{p_1}(\Omega,\mu)$, i.e., $$\norm{f}_p\leq \norm{f}_{p_0}^{1-\theta}\norm{f}_{p_1}^\theta\leq \norm{f}_{L^{p_0}\cap L^{p_1}}:=\operatorname{max}\{\norm{f}_{p_0},\norm{f}_{p_1}\}.$$ In particular, it implies that for $p_0=1$ and $p_1=\infty$, $L^1(\Omega,\mu)\cap L^\infty(\Omega,\mu)\xhookrightarrow{}L^p(\Omega,\mu)$. This behavior of the Lebesgue spaces with respect to the couple $\left(L^1(\Omega,\mu),L^\infty(\Omega,\mu)\right)$ and the  classical theorems of Marcinkiewicz and Riesz-Thorin \cite[Theorems~1.3.2 and 1.3.4]{Grafakos2014} motivate the precise definitions and main goals of the abstract interpolation theory, which can be summarized in the following diagram
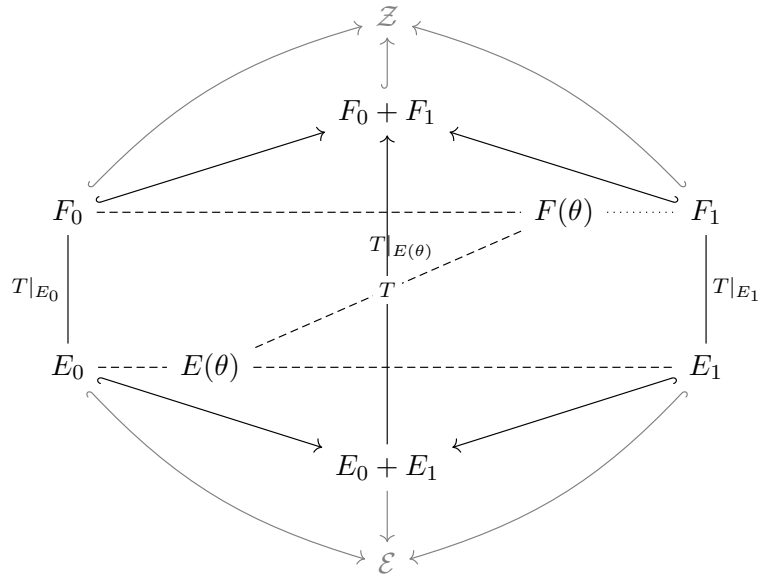
\[
\begin{tikzcd}
	&& {\color{gray} \mathcal{Z}} \\
	&& {F_0 + F_1} \\
	{F_0} &&& {F(\theta)} & {F_1} \\
	\\
	{E_0} & {E(\theta)} &&& {E_1} \\
	&& {E_0 + E_1} \\
	&& {\color{gray} \mathcal{E}}
	\arrow[hook, gray, from=2-3, to=1-3] 
	\arrow[bend left=15, hook, gray, from=3-1, to=1-3] 
	\arrow[hook, from=3-1, to=2-3]
	\arrow[dashed, no head, from=3-1, to=3-4]
	\arrow["{T|_{E_0}}"', -|>, no head, from=3-1, to=5-1]
	\arrow[dotted, no head, from=3-4, to=3-5]
	\arrow["{T|_{E(\theta)}}"'{pos=0.3}, dashed, -|>, no head, from=3-4, to=5-2]
	\arrow[bend left=-15, hook', gray, from=3-5, to=1-3] 
	\arrow[hook', from=3-5, to=2-3]
	\arrow["{T|_{E_1}}", -|>, no head, from=3-5, to=5-5]
	\arrow[dashed, no head, from=5-1, to=5-2]
	\arrow[hook', from=5-1, to=6-3] 
	\arrow[bend left=-15, hook', gray, from=5-1, to=7-3] 
	\arrow[dashed, no head, from=5-2, to=5-5]
	\arrow[hook, from=5-5, to=6-3] 
	\arrow[bend right=-15, hook, gray, from=5-5, to=7-3] 
	\arrow["T"{description}, from=6-3, to=2-3] 
	\arrow[gray, from=6-3, to=7-3] 
\end{tikzcd}
\]
Suppose that we have two pairs of normed spaces, which we will suppose complete for simplicity, $(E_0,E_1)$ and $(F_0,F_1)$, such that they are compatible in the sense of their sum and intersection being well defined as Banach spaces, which happens if both spaces of the couple are continuously embedded in a Hausdorff topological vector space. Such couples of Banach spaces are called \textit{compatible couples}. In our case $$E_0,E_1\xhookrightarrow{}\mathcal{E},\,F_0,F_1\xhookrightarrow{}\mathcal{Z}.$$ 
Let $T:E_0+E_1\to F_0+F_1$ be a bounded linear operator such that its restriction to $E_j$ is a bounded linear operator from $E_j\to F_j$, $j=0,1$. Such operators are called \textit{admissible operators} from $(E_0,E_1)\to (F_0,F_1)$. Suppose that we have two spaces $E(\theta), F(\theta)$ that are \textit{intermediate} with respect to the couples $(E_0,E_1), (F_0,F_1)$, in the sense that 
$$E_0\cap E_1\xhookrightarrow{}E(\theta)\xhookrightarrow{}E_0+E_1,\quad F_0\cap F_1\xhookrightarrow{}F(\theta)\xhookrightarrow{}F_0+F_1.$$ 
If we see the intermediate spaces between a compatible couple as a scale with endpoints being the members of such a couple, the parameter $\theta\in (0,1)$ is related to how far from the endpoints the intermediate space is. The main goal of interpolation theory is to construct intermediate spaces, and that the restriction of any admissible operator $T$ to $E(\theta)$ is again a bounded linear operator from $E(\theta)\to F(\theta)$, in such a way that its norm is related to the norm of the operator acting on the endpoints and the parameter $\theta$. These spaces, known as \textit{interpolation spaces}, are central to the theory. They retain many properties of the compatible couples (the \textit{endpoints} of the interpolation). Consequently, interpolation theory serves as a powerful tool for examining function spaces derived from the interpolation of compatible couples: which properties of compatible couples are inherited by interpolation spaces? That is the case of $L^p$ spaces with respect to couple $(L^1,L^\infty)$. For instance, for $\R^n$ equipped with the Lebesgue measure, and for every $1\leq p,q\leq \infty$, the couple $\left( L^p(\R^n),L^q(\R^n)\right)$ is compatible since we can choose $\mathcal{E}$ to be $L^0(\R^n)$, the space of equivalence classes of measurable functions with the topology generated by the metric 
$$d(f,g):=\int_{\R^n}\frac{|f(x)-g(x)|}{1+|f(x)-g(x)|}e^{-x^2}\,dx.$$ 

Given two compatible couples, various approaches, or \textit{interpolation methods}, can be utilized to derive interpolation spaces. These methods integrate seamlessly with the language of category theory, functioning as functors. In this context, we define the category of compatible couples of Banach spaces $\Cc\Cc $ as the one whose objects are compatible couples of Banach spaces and whose morphisms are the admissible operators between them. Then, an interpolation method will be a functor $\mathcal{F}:\Cc\Cc\to \textbf{Ban}$, where $\textbf{Ban}$ is the category of Banach spaces, such that given $(E_0,E_1),(F_0,F_1)\in \Cc\Cc$ and $T:(E_0,E_1)\to (F_0,F_1)$ admissible, we have $\mathcal{F}(T)=T|_{\mathcal{F}\left((E_0,E_1)\right)}:\mathcal{F}\left((E_0,E_1)\right)\to \mathcal{F}\left((F_0,F_1)\right)$ is linear and bounded. If we denote $E=\mathcal{F}\left((E_0,E_1)\right)$ and $F=\mathcal{F}\left((F_0,F_1)\right)$, the interpolation functors are classified in the following way: we say that $\mathcal{F}$ is \begin{itemize}
    \item \textit{exact}, if $$||T||_{\mathcal{L}(E,F)}\leq \operatorname{max}\left(||T||_{\mathcal{L}(E_0,F_0)},||T||_{\mathcal{L}(E_1,F_1)}\right),$$
    \item \textit{uniform}, if there exists a constant $C>1$ such that $$||T||_{\mathcal{L}(E,F)}\leq C\operatorname{max}\left(||T||_{\mathcal{L}(E_0,F_0)},||T||_{\mathcal{L}(E_1,F_1)}\right),$$
    \item \textit{exact of exponent $\theta\in (0,1)$}, if $$||T||_{\mathcal{L}(E,F)}\leq||T||_{\mathcal{L}(E_0,F_0)}^{1-\theta}||T||_{\mathcal{L}(E_1,F_1)}^\theta,$$
    \item \textit{of exponent $\theta\in (0,1)$}, if there exists a constant $C>1$ such that $$||T||_{\mathcal{L}(E,F)}\leq C||T||_{\mathcal{L}(E_0,F_0)}^{1-\theta}||T||_{\mathcal{L}(E_1,F_1)}^\theta.$$
\end{itemize} 

Another important notion of category theory with major applications in interpolation theory is the one of \textit{retractions} and \textit{coretractions}. Let $\Cc$ be some category and $A,B\in \operatorname{Obj}_{\Cc}$, where $\operatorname{Obj}_{\Cc}$ are the objects in the category $\Cc$. We say that $B$ is a \textit{retract} of $A$ if there exist $f\in \mor_{\Cc}(A,B)$ and $g\in \mor_{\Cc}(B,A)$ such that $f\circ g=i_B$, where $\mor_{\Cc}(A,B)$ stands for morphisms from $A$ to $B$. The morphism $f$ is called a \textit{retraction} and $g$ is the corresponding \textit{coretraction} or \textit{section}. In the category of Banach spaces $\textbf{Ban}$, given two Banach spaces $E,F$, the fact $F$ being a retract of $E$ means that there exist two bounded linear operators $R\in \mathcal{L}(E,F)$ and $S\in \mathcal{L}(F,E)$ such that $R\circ S$ is the identity operator. In the category of compatible couples of Banach spaces $\Cc\Cc$, a compatible couple $\overline{F}=(F_0,F_1)$ is a retract of another compatible couple $\overline{E}=(E_0,E_1)$ if there exist two admissible operators $R:\overline{E}\to \overline{F}$, $S:\overline{F}\to\overline{E}$, such that $R\circ S$ is the identity operator for $F_0+F_1$ with its restriction to $F_j$ being as well the identity on $\mathcal{L}(F_j)$, $j=0,1$, i.e., $$R|_{E_j}\circ S|_{F_j}=i_{F_j}.$$
Retracts in $\Cc\Cc$ are kept by interpolation functors in the following sense: if $\mathcal{F}$ is an interpolation functor and $\overline{F}$ is a retract of $\overline{E}$ in the sense of the category $\Cc\Cc$, then by the properties of $\mathcal{F}$, $\mathcal{F}(\overline{F})$ is a retract of $\mathcal{F}(\overline{E})$ in the sense of $\textbf{Ban}$. Moreover, retraction and coretraction are just the restrictions $R|_{\mathcal{F}(\overline{E})}$ and $S|_{\mathcal{F}(\overline{F})}$, respectively. Note that given $x\in \mathcal{F}(\overline{E})$, 
$$(S\circ R)^2 x=(S\circ R)\left((S\circ R)x\right)=S\circ (R\circ S)(Rx)=S(Rx)=(S\circ R)x,$$ 
hence $(S\circ R)^2=(S\circ R)$ is a projection in $\mathcal{F}(\overline{E})$. Since it is a continuous projection, its image is closed. We will denote by $E$ the image of the restriction of $S\circ R$ to $\mathcal{F}(\overline{E})$. For every $y\in \mathcal{F}(\overline{F})$, we have that $$Sy=S\left((R\circ S)y\right)=(S\circ R)(Sy),$$ so $Sy\in E$, and hence the operator $S:\mathcal{F}(\overline{F})\to E$ is well defined. Moreover, since any element $x\in E$ is identified by the image by the projection $S\circ R$ of $x\in \mathcal{F}(\overline{E})$, we have $x=(S\circ R)x=S(Rx)$, so the mapping is surjective. Finally, since $R\circ S$ is the identity for $\mathcal{F}(\overline{F})$, the mapping $S$ is injective. Then, according to the closed graph Theorem, $S:\mathcal{F}(\overline{F})\to E$ is an isomorphism, that is, $\mathcal{F}(\overline{F})$ is isomorphic to a complemented subspace of $\mathcal{F}(\overline{E})$. We have proved the following theorem.
\begin{Teor}[Retraction Theorem]\label{retract}
    Let $\overline{E},\overline{F}$ two compatible couples of Banach spaces such that $\overline{F}$ is a retract of $\overline{E}$. Let $R\in \mathcal{L}(\overline{E},\overline{F})$ be a retraction with corresponding coretraction $S\in \mathcal{L}(\overline{F},\overline{E})$. Then, given an interpolation functor $\mathcal{F}$, the interpolation space $\mathcal{F}(\overline{F})$ is isomorphic to a complemented subspace of the interpolation space $\mathcal{F}(\overline{E})$. The subspace is the image of the restriction of $S\circ R$ to $\mathcal{F}(\overline{E})$.
\end{Teor}
This theorem will be crucial in the computation of interpolation spaces, as it provides a method to transfer information from known interpolation spaces to others. Suppose that the interpolation space $\mathcal{F}(\overline{E})$ is well known for some particular interpolation functor $\mathcal{F}$ and compatible couple $\overline{E}$. If we are given another compatible couple $\overline{F}$ such that $\overline{F}$ is a retract, the construction of the operator $S$ will allow us to determine the interpolation space $\mathcal{F}(\overline{F})$ in the light of {Theorem \ref{retract}}, so the study of properties of the space $\mathcal{F}(\overline{F})$ is reduced to the study of retractions from the known space $\mathcal{F}(\overline{E})$.\\

For a comprehensive exposition of categorical interpolation theory, we refer to \cite{Brudnyi1991, Castillo2010}. For our purposes, we will focus on the real and complex interpolation methods and their interrelationships. 

\subsection{The Real Method of Interpolation}
There are several {\it real} methods of interpolation, many of them equivalent. The most popular is the one based on Peetre's $K$-functional, defined for a compatible couple of Banach spaces $(E_0,E_1)$ as $$K(t,x)=K\left(t,x;(E_0,E_1)\right):=\operatorname{inf}_{x=x_0+x_1}\left(\norm{x_0}_{E_0}+t\norm{x_1}_{E_1}\right),$$ for $x\in E_0+E_1$ and $t>0$. Clearly, it defines an equivalent norm on $E_0+E_1$ for every $t>0$. Now, given $1\leq p<\infty$ and $\theta\in (0,1)$, or $p=\infty$ and $\theta\in [0,1]$, we define the real interpolation spaces $(E_0,E_1)_{\theta,p}$ as the spaces 
$$(E_0,E_1)_{\theta,p}:=\{x\in E_0+E_1: \norm{x}_{\theta,p}<\infty\},$$ 
where 
$$\norm{x}_{\theta,p}:=\begin{cases}
    \left(\displaystyle\int_0^\infty \left(t^{-\theta}K(t,x)\right)^p\frac{dt}{t}\right)^{1/p},& 1\leq p<\infty,\\
    \operatorname{sup}_{t>0}t^{-\theta}K(t,x),& p=+\infty.
\end{cases}.$$
The functor $K_{\theta,p}:(E_0,E_1) \to (E_0,E_1)_{\theta,p}$ is an exact interpolation functor (\cite[Theorem~3.1.2]{BerghLofstrom1976} or \cite[Theorem~III.1.3]{GarciaSaez2024}), i.e., given an admissible operator ({\it morphism}) $T:(E_0,E_1)\to (F_0,F_1)$ between two compatible couples, such that 
$$\norm{T}_{\mathcal{L}(E_j,F_j)}=M_j,\,j=0,1,$$ 
then $T:(E_0,E_1)_{\theta,q}\to (F_0,F_1)_{\theta,q}$ is linear and bounded, with 
$$\norm{T}_{\mathcal{L}\left((E_0,E_1)_{\theta,q},(F_0,F_1)_{\theta,q}\right)}\leq M_0^{1-\theta}M_1^\theta.$$ 
A relevant consequence is that if we have $E_j\xhookrightarrow{}F_j$, $j=0,1$, by interpolating the embedding operator we get that $(E_0,E_1)_{\theta,q}\xhookrightarrow{}(F_0,F_1)_{\theta,q}$. Now, we enumerate the most important properties of real interpolation spaces.
\begin{Teor}[Properties of real interpolation spaces]\label{Kprops}
     Let $1\leq p\leq \infty$, $\theta\in (0,1)$ and $(E_0,E_1)$ be a compatible couple of Banach spaces. The space $(E_0,E_1)_{\theta,q}$ satisfies the following properties:
    \begin{itemize}
        \item[1.-] If $E_0,E_1$ are Banach spaces then $(E_0,E_1)_{\theta,p}$ is a Banach space.
        \item[2.-] $(E_0,E_1)_{\theta,p}=(E_1,E_0)_{1-\theta,p}$.
        \item[3.-] If $1\leq p_1\leq p_2\leq \infty$ then $$(E_0,E_1)_{\theta,1}\xhookrightarrow{}(E_0,E_1)_{\theta,p_1}\xhookrightarrow{}(E_0,E_1)_{\theta,p_2}\xhookrightarrow{}(E_0,E_1)_{\theta,\infty}.$$
        \item[4.-] If $E_0\xhookrightarrow{}E_1$ and $0<\theta_0<\theta_1<1$, then $$(E_0,E_1)_{\theta_0,p}\xhookrightarrow{}(E_0,E_1)_{\theta_1,p}.$$
        \item[5.-] Let $p<\infty$. If $E_0\cap E_1$ is dense in both $E_0$ and $E_1$ (that is called a regular couple), then it follows $$\left((E_0,E_1)_{\theta,p}\right)^*=(E_0^*,E_1^*)_{\theta,(1-1/p)^{-1}},$$ with equivalence of the norms. Moreover, if $E_0$ and $E_1$ are reflexive and $E_0\xhookrightarrow{}E_1$ with dense inclusion, then the space $(E_0,E_1)_{\theta,p}$ is reflexive.
    \end{itemize}
\end{Teor}
Detailed proofs of these facts can be found in \cite[Theorem~3.4.1, Theorem~3.4.2, Theorem~3.7.1]{BerghLofstrom1976} and \cite[Theorem~III.1.4, Corollary~III.1.6, Theorem~III.3.4, Theorem~III.3.9]{GarciaSaez2024}.\\

A major result in interpolation theory, and in particular for the real method, is the Reiteration Theorem, which implies that the real method is stable under reiteration of the real interpolation functor. However, the theorem applies for a larger class of spaces, the so-called $\theta$-class. Given $\theta\in[0,1]$ and a compatible couple $(E_0,E_1)$, we say that a Banach space $E$ is of \textit{$\theta$-class} for the couple if $E$ is an intermediate space and if we have the embeddings 
$$(E_0,E_1)_{\theta,1}\xhookrightarrow{}E\xhookrightarrow{}(E_0,E_1)_{\theta,\infty},\,\theta\in (0,1)$$ or $$\overline{E_\theta}^\Delta\xhookrightarrow{}E\xhookrightarrow{}(E_0,E_1)_{\theta,\infty},\,\theta\in \{0,1\},$$ where $\overline{E}^\Delta$ is the completion of $E_0\cap E_1$ under the norm of $E$. The embedding $(E_0,E_1)_{\theta,1}\xhookrightarrow{}E\xhookrightarrow{}(E_0,E_1)_{\theta,\infty}$ is equivalent to the existence of positive constants $C,c>0$ such that for every $x\in E_0\cap E_1$, $$ct^{-\theta}K(t,x)\leq \norm{x}_{E}\leq C\norm{x}_{E_0}^{1-\theta}\norm{x}_{E_1}^\theta,$$ which could be seen as the generalization of H\"older's inequality for intermediate spaces \cite[Proposition~III.5.3]{GarciaSaez2024}. Obviously, the endpoints $E_0$ and $E_1$ are respectively of class $0$ and $1$. Also, it is straightforward to prove that every interpolation functor of exponent $\theta\in (0,1)$ produces spaces in the $\theta$-class \cite[Lemma~III.6.16]{GarciaSaez2024}. The main result is the following (\cite[Theorem~3.11.5]{BerghLofstrom1976} or \cite[Theorem~III.5.6]{GarciaSaez2024}):  
\begin{Teor}[Reiteration theorem for the real method]\label{RealReiter}
    Let $(E_0,E_1)$ and $(X_0,X_1)$ be two compatible couples of Banach spaces such that $X_j$ is intermediate with respect to the couple $(E_0,E_1)$, $j=0,1.$ Given $0\leq \theta_0<\theta_1\leq 1$, $\alpha\in (0,1)$, $\theta:=(1-\alpha)\theta_0+\alpha\theta_1$, and $1\leq p\leq \infty$, it holds that if $X_i$ is of $\theta_i$-class, then $$(X_0,X_1)_{\alpha,p}=(E_0,E_1)_{\theta,p},$$ with equivalence of the norms.
\end{Teor}
We now present two important examples of spaces obtained by the real interpolation method exposed in this section that we shall require in the next section. 
\begin{Ejems}\label{ejemplosrealinterpol}
\end{Ejems}
\begin{enumerate}
        \item[i)] Let $E$ be a Banach space and $(U,\mu)$ a measure space with $\sigma$-finite positive measure. For $1\le p\le +\infty$, $L^p((U,\mu),E)$ stands for the Bochner space of $\mu$-measurable functions $f:U\to E$ such that
        \[ \|f\|_{L^p(U,E)}=\left( \int_U \|f(x)\|_E^p\,d\mu(x)\right)^\frac{1}{p}<+\infty.\]
        We denote this space by $L^p(E)$ is there is no risk of confusion. 
        
        For $f\in L^1(E)+L^\infty(E)$ we define $$\lambda(f,w):=\mu\left(\{t\in U: \norm{f(t)}_E>w\}\right),\,0<w<\infty,$$ the distribution function of $f$. We also define the decreasing rearrangement of $f$ as $$f^*(t):=\operatorname{inf}_{\lambda(f,w)\leq t} w,\,0<t<\infty.$$
        For $1\leq p<\infty$ and $1\leq q<\infty$
        we define the Lorentz space $$L^{p,q}(E)=\Bigg\{f\in L^1(E)+L^\infty(E):  \norm{f}_{L^{p,q}}:=\left(\int_0^\infty \left(t^{1/p}f^*(t)\right)^q\,\frac{dt}{t}\right)^{1/q}<\infty\Bigg\},$$ and for $q=\infty$, $$L^{p,\infty}(E):=\Bigg\{f\in L^1(E)+L^\infty(E): \norm{f}_{L^{p,\infty}}:=\operatorname{sup}_{t>0}t^{1/p}f^*(t)<\infty.\Bigg\}.$$ When $p=\infty=q$ it is conventional to set $L^{\infty,\infty}(E)=L^\infty(E)$. In general, those spaces are quasi-normed spaces, but for $p>1$ it is possible to replace the quasi-norm with a norm, which makes them complete. From the definition, it is straightforward to see that $L^{p,p}(E)=L^p(E)$ with equality of norms for $1\leq p\leq \infty$. Lorentz spaces can be obtained by real interpolation of Lebesgue spaces. In fact, for $0<\theta<1$ and $1\leq q\leq \infty$, 
        $$\left(L^1(E),L^\infty(E)\right)_{\theta,q}=L^{\frac{1}{1-\theta},q}(E),$$ 
        see \cite[Theorem~1.18.6.1]{Triebel1992} for a proof. By {Theorem \ref{RealReiter}} it follows that for $\theta\in (0,1)$, $1<p_0<p_1<\infty$, and $1\leq q,\,q_0,\,q_1\leq \infty$, 
        $$\left(L^{p_0,q_0}(E),L^{p_1,q_1}(E)\right)_{\theta,q}=L^{p(\theta),q}(E),\,\frac{1}{p(\theta)}=\frac{1-\theta}{p_0}+\frac{\theta}{p_1}.$$ 
        In particular, $$\left(L^{p_0}(E),L^{p_1}(E)\right)_{\theta,p(\theta)}=L^{p(\theta)}(E).$$ Moreover, given a compatible couple of Banach spaces $(E_0,E_1)$, $$\left(L^{p_0}(E_0),L^{p_1}(E_1)\right)_{\theta,p(\theta)}=L^{p(\theta)}\left((E_0,E_1)_{\theta,p(\theta)}\right),$$ see \cite[Theorem~5.6.2]{BerghLofstrom1976}. Interestingly, there is no generalization for $q\not =p(\theta)$ and in general, we only have the embedding $$\left(L^{p_0}(E_0),L^{p_1}(E_1)\right)_{\theta,q}\xhookrightarrow{}L^{p(\theta)}\left((E_0,E_1)_{\theta,q}\right),$$
        (see \cite{Cwikel1974}).
\item[ii)] Let $s\in (0,1)$ and $p\in [1,\infty)$. We define the Gagliardo space $W^{s,p}(\R^n)$ as 
$$W^{s,p}(\R^n):=\Bigg\{u\in L^p(\R^n): \frac{|u(x)-u(y)|}{|x-y|^{\frac{n}{p}+s}}\in L^p(\R^n\times \R^n)\Bigg\},$$ 
with the norm 
$$\norm{u}_{W^{s,p}}:=\norm{u}_p+[u]_{W^{s,p}},$$ 
where $$[u]_{W^{s,p}}:=\left(\int_{\R^n}\int_{\R^n}\frac{|u(x)-u(y)|^p}{|x-y|^{n+sp}}\,dx\,dy\right)^\frac{1}{p},$$ is the Gagliardo seminorm. In the case $p=\infty$, we have the natural identification of  $W^{s,\infty}(\R^n)$ with $C^{0,s}(\R^n)$.

Gagliardo spaces are obtained by means of real interpolation of the couple $\left(L^p(\R^n),W^{1,p}(\R^n)\right)$. In particular we have the equality $$W^{s,p}(\R^n)=\left(L^p(\R^n),W^{1,p}(\R^n)\right)_{s,p},\,s\in (0,1),\,1<p<\infty,$$ with equivalence of the norms. We refer to \cite[Theorem~12.5]{Leoni2023} and \cite[Example~1.8]{Lunardi2018} for detailed proofs of these facts, and to \cite{DiNezzaPalatucciValdinoci2012} to a complete discussion of Gagliardo spaces without interpolation theory.
\end{enumerate}

\subsection{The Complex Method of Interpolation}
The other major method of interpolation is the Complex Method of Interpolation. Briefly introduced by Lions \cite{Lions1960} and Calder\'on \cite{Calderon1963}, it was later developed in depth in the seminal paper of the second one \cite{Calderon1964}. The ideas of Calder\'on have been generalized by several authors, see for example the work of Schechter \cite{Schechter1967}. The complex interpolation method is based on the theory of holomorphic vector valued functions, in particular the ones defined on the strip $$S=\{z\in\C: 0\leq\operatorname{Re}z\leq 1\},$$ taking values on some complex Banach space $E$. From the Hadamard's three line Theorem and the maximum modulus principle it follows that for any function $f:S\to E$ holomorphic on $\operatorname{int}S$ and  continuous and bounded in $S$, the inequality 
$$\norm{f(z)}_E\leq \operatorname{max}\{\operatorname{sup}_{t\in\R}\norm{f(it)}_E, \operatorname{sup}_{t\in\R}\norm{f(1+it)}_E\},\,z\in S$$ 
holds. This motivates the following definition:\\

Given a compatible couple of Banach spaces $(E_0,E_1)$, we define the space $\mathfrak{F}\left((E_0,E_1)\right)$ as the space of functions $f:S\to E_0+E_1$ such that the following hold:
\begin{itemize}
    \item $f$ is holomorphic in $\operatorname{int}S$, and continuous and bounded in $S$.
    \item The functions $t\mapsto f(j+it)$ are continuous from $\R\to E_j$, $j=0,1$, and such that $$\norm{f(j+it)}_{E_j}\to 0,\,j=0,1,$$ as $|t|\to \infty$.
\end{itemize}
Clearly, $\mathfrak{F}(E_0,E_1)$ is a vector space. Moreover, endowed with the norm
$$\norm{f}_{\mathfrak{F}(\overline{E})}:=\operatorname{max}\{\operatorname{sup}_{t\in\R}\norm{f(it)}_{E_0},\operatorname{sup}_{t\in\R}\norm{f(1+it)}_{E_1}\},\,f\in \mathfrak{F}(\overline{E}),$$ it becomes a Banach space. Then, for any $\theta\in [0,1]$, we define $[E_0,E_1]_\theta$ as the space of all $x\in E_0+E_1$ such that there exists $f\in \mathfrak{F}(\overline{E})$ with $f(\theta)=x$, and with the norm 
$$\norm{x}_{\theta}:=\operatorname{inf}\{\norm{f}_{\mathfrak{F}(\overline{E})}: f\in \mathfrak{F}(\overline{E}), f(\theta)=x\}.$$ 
The functor, 
$$\Cc_\theta:(E_0,E_1)\to [E_0,E_1]_\theta,$$ is an exact interpolation functor of exponent $\theta$, i.e., the space $[E_0,E_1]_\theta$ is a Banach space intermediate with respect to the couple $(E_0,E_1)$, and for any admissible operator $T:(E_0,E_1)\xhookrightarrow{}(F_0,F_1)$ between two compatible couples, such that 
$$\norm{T}_{\mathcal{L}(E_j,F_j)}=M_j,\,j=0,1,$$ then $T:[E_0,E_1]_{\theta}\to [F_0,F_1]_{\theta}$ with 
$$\norm{T}_{\mathcal{L}\left([E_0,E_1]_{\theta},[F_0,F_1]_{\theta}\right)}\leq M_0^{1-\theta}M_1^\theta.$$ 
A relevant consequence is that if we have $E_j\xhookrightarrow{}F_j$, $j=0,1,$ interpolating the embedding operator, we get that $[E_0,E_1]_{\theta}\xhookrightarrow{}[F_0,F_1]_{\theta}$. Another interesting consequence is that the spaces $[E_0,E_1]_\theta$ are of $\theta$-class, and hence, for $\theta_0,\theta_1\in [0,1]$, $\alpha\in (0,1)$ and $1\leq p\leq \infty$, {Theorem \ref{RealReiter}} holds, hence 
$$\left([E_0,E_1]_{\theta_0},[E_0,E_1]_{\theta_1}\right)_{\alpha,p}=(E_0,E_1)_{\theta,p},$$ where $\theta=(1-\alpha)\theta_0+\alpha\theta_1$. \\

Now we enumerate the most important properties of such spaces.
\begin{Teor}[Properties of Complex interpolation spaces]\label{ComplexProps}
        Let $(E_0,E_1)$ a compatible couple of Banach spaces, and $\theta\in [0,1]$. Then, we have
        \begin{enumerate}
            \item[1.-] $[E_0,E_1]_{\theta}=[E_1,E_0]_{1-\theta}$.
            \item[2.-]If $E_0\xhookrightarrow{}E_1$ and $\theta_0<\theta_1$, $[E_0,E_1]_{\theta_0}\xhookrightarrow{}[E_0,E_1]_{\theta_1}$.
            \item[3.-] If $E_0=E_1$ and $0<\theta<1$, $[E_0,E_1]_\theta=E_0$.
            \item[4.-] $E_0\cap E_1$ is dense in $[E_0,E_1]_\theta$.
            \item[5.-] $[E_0,E_1]_j$ is a closed subspace of $E_j$ with coincidence of the norm in $[E_0,E_1]_j$, $j=0,1$.
            \item[6.-] If $(E_0,E_1)$ is a regular couple and at least one of $E_0$ or $E_1$ is reflexive, then $$\left([E_0,E_1]_\theta\right)^*=[E_0^*,E_1^*]_\theta,\,\theta\in (0,1).$$
            \item[7.-]  If at least one of $E_0$ or $E_1$ is reflexive, then the space $[E_0,E_1]_\theta,\,\theta\in (0,1),$ is reflexive
         \end{enumerate}
\end{Teor}
Detailed proofs of these facts can be found in \cite[Theorem~4.2.1, Theorem~4.2.2, Theorem~4.5.1]{BerghLofstrom1976} and \cite[Proposition~IV.1.8, Theorem~IV.5.4, Theorem~IV.5.6]{GarciaSaez2024}.
There are also two more reiteration results for the complex method, concerning the complex interpolation of real interpolation spaces and the reiteration of complex spaces under complex interpolation. The first result is due to Karadzov \cite{Karadzov1974}, while the second was proved by Cwikel in \cite{Cwikel1978}, improving the original result by Calder\'on.
\begin{Teor}\label{ComplexReiter}
    Let $(E_0,E_1)$ be a compatible couple of Banach spaces, $\theta_0,\theta_1\in (0,1)$, $\alpha\in (0,1)$ and $1\leq p_0,p_1\leq \infty$. Define
     $$\theta(\alpha)=(1-\alpha)\theta_0+\alpha\theta_1,\,\frac{1}{p(\alpha)}=\frac{1-\alpha}{p_0}+\frac{\alpha}{p_1}.$$ Then, \begin{itemize}
        \item $\left[(E_0,E_1)_{\theta_0,p_0},(E_0,E_1)_{\theta_1,p_1}\right]_\alpha=(E_0,E_1)_{\theta(\alpha),p(\alpha)}.$
        \item $\left[[E_0,E_1]_{\theta_0},[E_0,E_1]_{\theta_1}\right]_\alpha=[E_0,E_1]_{\theta(\alpha)}.$
    \end{itemize}
\end{Teor}
As we did for the real method, we present some useful examples of spaces obtained by complex interpolation.
\begin{Ejems}\label{ComplexEjemes}
\end{Ejems}
\begin{itemize}
    \item[i)] The following interpolation identity may be understood as the generalization of the Riesz-Thorin theorem. Let $E$ a Banach space, $1\leq p_0<\infty$ and $1\leq p_1\leq \infty$, and $\theta\in (0,1)$. Then, $$[L^{p_0}(E),L^{p_1}(E)]_\theta=L^{p(\theta)}(E),\,\frac{1}{p(\theta)}=\frac{1-\theta}{p_0}+\frac{\theta}{p_1},$$ see \cite[Theorem~1.18.6.2]{Triebel1995}. Moreover, given a compatible couple of Banach spaces $(E_0,E_1)$, for $1\leq p_0,p_1<\infty$,     $$[L^{p_0}(E_0),L^{p_1}(E_1)]_\theta=L^{p(\theta)}\left([E_0,E_1]_\theta\right),$$ see \cite[Theorem~5.1.2]{BerghLofstrom1976}.
    
    \item[ii)] Let $s\in \R$, $0<q\leq \infty$ and $E$ be a Banach space. We define $\ell^{s,q}(E)$ as the space of $E$-valued sequences $(x_m)_{m=0}^\infty$ such that 
    \begin{align*}\norm{(x_m)}_{\ell^{s,q}(E)}&:=\left(\sum_{m=0}^\infty \left(2^{ms}\norm{x_m}_E\right)^q\right)^{1/q},\quad q<\infty,\\
    \norm{(x_m)}_{\ell^{s,\infty}(E)}&:=\sup_{m\geq 0} 2^{ms}\norm{x_m}_E,\quad q=\infty.\end{align*} 
     Then, given $s_0,s_1\in \R$, $1\leq q_0,q_1\leq \infty$, $\theta\in (0,1)$ and $(E_0,E_1)$ a compatible couple of Banach spaces, we have that $$[\ell^{s_0,q_0}(E_0),\ell^{s_1,q_1}(E_1)]_\theta=\ell^{s(\theta),q(\theta)}\left([E_0,E_1]_\theta\right),$$
    where,
    $$s(\theta)=(1-\theta)s_0+\theta s_1,\,\frac{1}{q(\theta)}=\frac{1-\theta}{q_0}+\frac{\theta}{q_1},$$ 
    see \cite[Theorem~5.6.3]{BerghLofstrom1976}.
    \item[iii)] Consider the space $\textbf{BMO}(\R^n)$ of functions with \textit{bounded mean oscillation}, (see \cite{JohnNirenberg1961}), i.e., the space of locally integrable functions $f$ on $\R^n$ such that $$\operatorname{sup}_Q\frac{1}{|Q|}\int_Q|f(x)-f_Q|\,dx=:\norm{f}_{\textbf{BMO}(\R^n)}<\infty,$$ where the supremum ranges over all cubes $Q\subset \R^n$, with $|Q|$ the Lebesgue measure of the cube and $f_Q=1/|Q|\int_Qf$. We will use this space in the sequel as an endpoint for interpolation results. The crucial thing is that $\textbf{BMO}(\R^n)$ behaves like $L^\infty(\R^n)$ for the generalized Riesz-Thorin theorem. This space is connected with the real variable Hardy space $\mathcal{H}^1(\R^n)$ (see \cite{FeffermanStein1972}), for which it is known that (see \cite{JansonJones1982} or \cite{FeffermanStein1972}) $$[\mathcal{H}^1(\R^n),L^p(\R^n)]_\theta=L^q(\R^n),\,\frac{1}{q}=(1-\theta)+\frac{\theta}{p},\,\theta\in(0,1).$$ Now, since $L^p(\R^n)$ is reflexive, {Theorem \ref{ComplexProps} (vi)} holds, and using the fact that $\left(\mathcal{H}^1(\R^n)\right)^*=\textbf{BMO}(\R^n)$ (\cite[Theorem~2]{FeffermanStein1972}), we have that \begin{align*}
    L^{q'}(\R^n)&=\left([\mathcal{H}^1(\R^n),L^{p}(\R^n)]_\theta\right)^*=\left[\left(\mathcal{H}^1(\R^n)\right)^*,\left(L^{p}(\R^n)\right)^*\right]_\theta\\
    &=[\textbf{BMO}(\R^n),L^{p'}(\R^n)]_\theta=[L^{p'}(\R^n),\textbf{BMO}(\R^n)]_{1-\theta},
\end{align*}
where $1/p+1/p'=1$ and $1/q+1/q'=1$. Taking $\theta=1-s$, we have that 
$$[L^{p'}(\R^n),\textbf{BMO}(\R^n)]_s=L^{q'}(\R^n),$$ 
where
$$\frac{1}{q'}=1-\frac{1}{q}=1-s-\frac{1-s}{p}=1-s-(1-s)\left(1-\frac{1}{p'}\right)=\frac{1-s}{p'}.$$ Then, we get that $$[L^p(\R^n),\textbf{BMO}(\R^n)]_s=L^{q(s)}(\R^n),$$ where $\frac{1}{q(s)}=\frac{1-s}{p}$.
\end{itemize}

\subsection{Relationship between both methods}
 Since the complex interpolation functor $\Cc_\theta$ is an exact interpolation functor of exponent $\theta$, it holds 
 $$(E_0,E_1)_{\theta,1}\xhookrightarrow{}\Cc_\theta(E_0,E_1)\xhookrightarrow{}(E_0,E_1)_{\theta,\infty},\quad 0<\theta<1,$$ for any compatible couple of Banach spaces $(E_0,E_1)$. A natural question arises: could any of the indices $1$ and $\infty$ be replaced by some $1<p<\infty$? In general, real and complex interpolation yield different results \cite[6.23]{BerghLofstrom1976}, but with some extra hypotheses over the spaces $E_0$ and $E_1$, we can show some nice embeddings between the two methods. This question was first addressed by Peetre \cite{Peetre1969}, where he introduced the following concept.

Given a Banach space $E$ and the Bochner space $L^p(\R,E)$, which we denote by $L^p(E)$, the Fourier transform on $L^p(E)$ is defined by
 $$\mathcal{F}\{f\}(\xi)=\int_{\R} f(x)e^{-2\pi ix\cdot \xi}\,dx,$$ as a Bochner integral.
 
 \begin{Defi}\em
 Let $E$ a Banach space and $1\leq p\leq 2$. We say that $E$ is of \textit{$p$-type} if the Fourier transform is a linear bounded operator from $L^p(E)$ to $L^q(E)$, where $1=1/p+1/q$.
 \end{Defi}

 Observe that every Banach space $E$ is of $1$-type, since the Fourier Transform of $f\in L^1(E)$ satisfies $\norm{\mathcal{F}f}_\infty<\infty$, so $\mathcal{F}\in \mathcal{L}\left(L^1(E),L^\infty(E)\right)$. It is also straightforward that every Hilbert space $H$ is of $2$-type due to Plancherel Identity. 

Moreover, all $L^p$ spaces are of some type. In particular, for $1<p<\infty$, $L^p(\Omega,d\mu)$ is of type $\operatorname{min}\{p,q\}$, $1=1/p+1/q$, where $(\Omega,\mu)$ is some measure space. Indeed, let $1<p\leq 2$ and consider the space $L^p\left(L^p(\Omega,d\mu)\right)$. For a given $f\in L^p\left(L^p(\Omega,d\mu)\right)$, $$\mathcal{
F}\{f\}(\xi)=\int_\R e^{-2\pi ix\xi}f(x)\,dx,$$ in the Bochner integral sense. Then, for every $\xi\in\R$, $\mathcal{F}\{f\}(\xi)$ is a function defined on $\Omega$ whose $p$-norm is finite. Using the Minkowski integral inequality \begin{align*}
    \left(\int_\R\norm{\mathcal{F}\{f\}(\xi)}_p^q\,d\xi\right)^{p/q}&=\left(\int_\R\left(\int_\Omega|\mathcal{F}\{f\}(\xi)(s)|^p\,d\mu(s)\right)^{q/p}\right)^{p/q}\\
    &\leq \int_\Omega\left(\int_\R |\mathcal{F}\{f\}(\xi)(y)|^q\,d\xi\right)^{p/q}d\mu(s), 
\end{align*}
and by the Hausdorff-Young theorem, $$\left(\int_\R|\mathcal{F}\{f\}(\xi)(s)|^q\,d\xi\right)^{1/q}\leq \left(\int_\R |f(x)(s)|^p\,dx\right)^{1/p},$$ thus \begin{align*}
    \left(\int_\R\norm{\mathcal{F}\{f\}(\xi)}_p^q\,d\xi\right)^{p/q}& \leq \int_\Omega d\mu(s)\int_\R |f(x)(s)|^p\,ds=\int_\R \norm{f(x)}_p^p\,dx,
\end{align*} so $$\norm{\mathcal{F}\{f\}}_{L^q\left(L^p(\Omega,d\mu)\right)}\leq \norm{f}_{L^p\left(L^p(\Omega,d\mu)\right)},$$ so $L^p(\Omega,d\mu)$ is of $p$-type. In the case $2\leq p<\infty$, with analogous reasoning we get that $L^p(\Omega,d\mu)$ is of $q$-type.\\

The previous concept relates to the interpolation theory due to the following result. We include its proof for readers' convenience. 
\begin{Lema}
    Let $(E_0,E_1)$ a compatible couple of Banach spaces such that $E_i$ is of $p_i$-type, $1\leq p_i\leq 2$, $i=0,1$. Then, for any $\theta\in (0,1)$, $(E_0,E_1)_{\theta,p(\theta)}$ and $[E_0,E_1]_\theta$ are of $p(\theta)$-type, where $$\frac{1}{p(\theta)}=\frac{1-\theta}{p_0}+\frac{\theta}{p_1}.$$
\end{Lema}
\noindent\textbf{Proof:}
We have that the Fourier transform $\mathcal{F}\in\mathcal{L}\left(L^{p_i}(E_i),L^{q_i}(E_i)\right)$ for $i=0,1$, hence 
$$\mathcal{F}\in \mathcal{L}\left(\left(L^{p_0}(E_0),L^{p_1}(E_1)\right)_{\theta,p(\theta)},\left(L^{q_0}(E_0),L^{q_1}(E_1)\right)_{\theta,p(\theta)}\right),$$ and $$\mathcal{F}\in \mathcal{L}\left([L^{p_0}(E_0),L^{p_1}(E_1)]_{\theta},[L^{q_0}(E_0),L^{q_1}(E_1)]_{\theta}\right)$$ 
By \cite[Theorem~5.6.2]{BerghLofstrom1976} $$\left(L^{p_0}(E_0),L^{p_1}(E_1)\right)_{\theta,p(\theta)}=L^{p(\theta)}\left((E_0,E_1)_{\theta,p(\theta)}\right),$$ and 
$$\left(L^{q_0}(E_0),L^{q_1}(E_1)\right)_{\theta,p(\theta)}\xhookrightarrow{}L^{q(\theta)}\left((E_0,E_1)_{\theta,p(\theta)}\right),$$ (see \cite{Cwikel1974}) where 
$$\frac{1}{q(\theta)}=\frac{1-\theta}{q_0}+\frac{\theta}{q_1},\mbox{  so  }1=\frac{1}{p(\theta)}+\frac{1}{q(\theta)}.$$ 
We conclude that the space $(E_0,E_1)_{\theta,p(\theta)}$ is of $p(\theta)$-type. Also, by \cite[Theorem~5.1.2]{BerghLofstrom1976}
$$[L^{p_0}(E_0),L^{p_1}(E_1)]_\theta=L^{p(\theta)}\left([E_0,E_1]_\theta\right)$$ 
and 
$$[L^{q_0}(E_0),L^{q_1}(E_1)]_{\theta}=L^{q(\theta)}\left([E_0,E_1]_\theta\right),$$ 
so we conclude that $[E_0,E_1]_\theta$ is of $p(\theta)$-type.\qed

With this notion, Peetre (\cite{Peetre1969}) showed the following major result connecting the Real and Complex interpolation spaces .
\begin{Teor}\label{PeetreTypeTheor}
    Let $(E_0,E_1)$ be a compatible couple of Banach spaces such that $E_i$ is of $p_i$-type, $1\leq p_i\leq 2$, $i=0,1$, with conjugate exponents $q_i$ ($1=1/p_i+1/q_i$), $i=0,1$. Then, $$(E_0,E_1)_{\theta,p(\theta)}\xhookrightarrow{}[E_0,E_1]_\theta\xhookrightarrow{}(E_0,E_1)_{\theta,q(\theta)},\,0<\theta<1,$$
    where
    $$\frac{1}{p(\theta)}=\frac{1-\theta}{p_0}+\frac{\theta}{p_1},\quad \frac{1}{q(\theta)}=\frac{1-\theta}{q_0}+\frac{\theta}{q_1}.$$
\end{Teor}
Peetre's original proof is challenging to follow. Therefore, for the readers' interest, we provide a comprehensive proof of this theorem in {Appendix \ref{B}}.

\begin{Obse}\em\label{CoincidenciaDeMetodos}
    Since every Banach space is of $1$-type, we directly obtain the fact that $$(E_0,E_1)_{\theta,1}\xhookrightarrow{}\Cc_\theta\overline{E}\xhookrightarrow{}(E_0,E_1)_{\theta,\infty},\,\theta\in (0,1),$$ for any compatible couple of Banach spaces.
    
     Another important consequence of the theorem is that the real and complex interpolation methods yields the same spaces with equivalent norms when we restrict ourselves to Hilbert spaces. In fact, given a compatible couple of Hilbert spaces $(X,Y)$, i.e, as normed spaces they form a compatible couple of Banach spaces, since $X$ and $Y$ are of $2$-type, and the conjugate of $2$ is $2$, we have that for any $\theta\in (0,1)$, 
     $$(X,Y)_{\theta,2}\xhookrightarrow{}[X,Y]_\theta\xhookrightarrow{}(X,Y)_{\theta,2},$$so 
     $$(X,Y)_{\theta,2}=[X,Y]_\theta.$$
\end{Obse}
\section{Bessel Potentials and Complex Interpolation of Sobolev spaces}

Building on the interpolation concepts discussed in the previous section, we now turn our attention to the complex interpolation of Sobolev spaces. This section focuses on the classical Bessel potential spaces, which are derived through interpolation using the Bessel potential. The primary result is a straightforward proof demonstrating the equivalence between Bessel spaces and the complex interpolation space $[L^p(\R^n),W^{k,p}(\R^n)]_\theta$. Leveraging the well-established properties of the $\Cc_\theta$ functor, we can effortlessly derive the main characteristics of these spaces.

\subsection{Complex interpolation and classical Bessel potential spaces}

The family of Bessel potential spaces is a cornerstone in Functional Analysis, particularly in the study of partial differential equations. These spaces are known by various names in the literature, with \textit{Bessel potential spaces} being the most widely used term. Initially introduced by Aronszajn and Smith \cite{Aronszajn1961} for the Hilbertian case $p=2$, these spaces have been extensively studied due to their significance in the theory of partial differential equations. Although Aronszajn and Smith introduced these spaces, they named them in honor of Friedrich Wilhelm Bessel, as the Bessel potential of a function is the convolution of that function with the Bessel kernel, which in turn can be represented using modified Bessel functions (see \cite{Aronszajn1961} for more details).

Notable works by Lions and Magenes \cite{LionsMagenes1972,LionsMagenes1961} explored interpolation methods, showing their equivalence to Calder\'on's method for the Hilbertian case. Calderón himself later extended these ideas to the range $1\le p\le \infty$ in \cite{Calderon1961}, defining them as the image of $L^p$ under the Bessel potential. Calderón's work is deeply rooted in Harmonic Analysis, and it was not until Lions' work in \cite{Lions1960} that these spaces were considered as interpolation spaces. However, Lions merely mentioned them as an example of the new complex interpolation method, without providing further proofs or commentaries, unlike Calderón in \cite{Calderon1963}.

Recently, Campos and Rodrigues \cite{Campos2023,campos2023Primero,Campos2024} have proposed naming these spaces \textit{Lions-Calderón spaces} in recognition of their complete introduction by them. They are also sometimes referred to as \textit{General Sobolev spaces} as noted in \cite{GrafakosII2014}. In Russian literature, these spaces are known as \textit{Liouville spaces}, as pointed out by Triebel in \cite{Triebel1992}. From the perspective of fractional gradients \cite{BellidoCuetoMoraCorral2021,BellidoCuetoMoraCorral2020,ComiStefani2019,BrueCalziComiStefani2022,ComiStefani2023,ShiehSpector2015}, the terms \textit{Fractional Sobolev spaces} and \textit{Distributional Fractional Sobolev spaces} are used. It is arguable that these terms are more appropriate for them than for the spaces $W^{s,p}$
derived from Gagliardo seminorms, as fractional gradients generalize the classical gradient \cite{BellidoCuetoMoraCorral2021}. Given the classical use of \textit{Fractional Sobolev spaces} for $W^{s,p}$
spaces and the infrequent use of other terms, we will maintain the classical denomination of Bessel potential spaces for our discussion.

\begin{Defi}\em[Bessel potential]
    Let $s\in \C$, $f\in \mathcal{S}'(\R^n)$ (where $ \mathcal{S}'(\R^n)$ is the space of tempered distributions) and $\xi\in \R^n$. We define the Bessel potential $\Lambda_s$ of order $s$ of $f$ as $$\Lambda_{s}f=\mathfrak{F}^{-1}\left((1+4\pi^2|\xi|^2)^{-s/2}\mathfrak{F}f(\xi)\right).$$
\end{Defi}
\begin{Obse}\em\label{BesselLpLp}
Note that $\Lambda_s$ is well-posed. Indeed, since $\left(1+4\pi^2|\xi|^2\right)^{-s/2}$ is a smooth function with all its derivatives with at most polynomial growth, its product with a tempered distribution defines another tempered distribution. Also, since the Fourier transform and its inverse transform tempered distributions into tempered distributions, 
$$\Lambda_s:\mathcal{S}'(\R^n)\to \mathcal{S}'(\R^n); \quad f\mapsto \mathfrak{F}^{-1}\left((1+4\pi^2|\xi|^2)^{-s/2}\mathfrak{F}f(\xi)\right).$$ 
is well defined as a linear operator. Moreover, from the basic properties of the Fourier transform acting on tempered distributions, we can easily derive that the Bessel potential has the semigroup property with the composition, i.e., for any $s_0,s_1\in \C$ and $f\in \mathcal{S}'(\R^n)$, 
$$\Lambda_{s_0+s_1}f=\Lambda_{s_1}(\Lambda_{s_0}f).$$ 

It can be shown that for every $s>0$,
\[\left(1+4\pi^2|\xi|^2\right)^{-s/2}= \mathcal{F}G_s, \]
with
$$G_s(x)=\frac{1}{(4\pi)^{n/2}\Gamma(n/2)}\int_0^\infty e^{-t/(4\pi)}e^{-|x|^2\pi/t}t^{(s-n)/2}\,\frac{dt}{t},\,x\in \R^n,$$ 
which satisfies that $$\norm{G_s}_1=\int_{\R^n}G_s(x)\,dx=1,$$ see \cite{Hao2016,Mizuta1996}. Consequently,
\[\Lambda_sf=G_s*f.\]
From there, we can derive a useful result: for every $s>0$ and $1\le p \le \infty$, $\Lambda_s$ is a continuous linear operator from $L^p(\R^n)$ into itself, with norm at most one. This follows from Young's convolution inequality, since $$\norm{\Lambda_sf}_p=\norm{G_s*f}_p\leq \norm{G_s}_1\norm{f}_p=\norm{f}_p,$$
therefore,
$$\norm{\Lambda_s}_{\mathcal{L}\left(L^p(\R^n),L^p(\R^n)\right)}\leq 1.$$
\end{Obse}
We now present the Bessel potential spaces as they were introduced by Lions and Calder\'on, respectively.
\begin{Defi}\em[Bessel potential space (Lions)]
    Let $s\in \R$ and $1\leq p\leq \infty$. We define the \textit{Bessel potential space (Lions)} $\Lambda_{s,p}(\R^n)$ as $$\Lambda^{s,p}(\R^n):=\{f \in \mathcal{S}'(\R^n): \Lambda_{-s}f\in L^p(\R^n)\},$$ with the norm $$\norm{f}_{\Lambda^{s,p}}:=\norm{\Lambda_{-s}f}_p.$$
\end{Defi}
Notice that $\Lambda_s$ stands for the Bessel potential, and $\Lambda^{s,p}$ for the Bessel potential space with the definition given by Lions. 

\begin{Defi}\em[Bessel Potential space (Calder\'on)] Let $s\in\R$ and $1\leq p\leq\infty$. We define the \textit{Bessel Potential space (Calder\'on)} $L_s^p(\R^n)$ as the image of $L^p(\R^n)$ by $\Lambda_{s}$, i.e., $$L_s^p(\R^n):=\{\Lambda_{s}f: f\in L^p(\R^n)\},$$ with the norm $$\norm{\Lambda_{s}f}_{L_s^p}=\norm{f}_p.$$
\end{Defi}
Obviously, for every $s\in \R$ and $1\leq p\leq \infty$, the spaces $\Lambda^{s,p}(\R^n)$ and $L_s^p(\R^n)$ are equal, with equality of norms.

One of the crucial things about Bessel potential spaces is that they generalize the classical Sobolev spaces, i.e., for any $1<p<\infty$ and $k\in \mathbb{N}$, we have that $$\Lambda^{k,p}(\R^n)=W^{k,p}(\R^n),$$ as proved by Calder\'on \cite[Theorem~7]{Calderon1961}. See also \cite[Chapter~7, Theorem~2.2]{Mizuta1996}.

From the definitions of Lions and Calder\'on, the following result of density is straightforward.
\begin{Prop}\label{SchwartzDenseBessel}
    Let $s>0$ and $1\le p<\infty$. Then, the space $\mathcal{S}(\R^n)$ is dense in $\Lambda^{s,p}(\R^n)$.
\end{Prop}
\noindent\textbf{Proof:} Let $f\in \Lambda^{s,p}(\R^n)$, which means that $\Lambda_{-s}f\in L^p(\R^n)$. Since $\mathcal{S}(\R^n)$ is dense in $L^p(\R^n)$, for every $\varepsilon>0$ there exists $g\in \mathcal{S}(\R^n)$ such that $$\norm{\Lambda_{-s}f-g}_p<\varepsilon.$$ Since $\Lambda_s$ is a bounded linear operator from $\mathcal{S}(\R^n)$ into itself, $$\norm{\Lambda_{-s}f-g}_p=\norm{f-\Lambda_s g}_{\Lambda^{s,p}},$$ and hence $$\norm{f-\Lambda_sg}_{\Lambda^{s,p}}<\varepsilon,$$ which means that $\mathcal{S}(\R^n)$ is dense in $\Lambda^{s,p}(\R^n)$. \qed

Another immediate but crucial fact is that $\Lambda_{-t}$ is an isomorphism between $\Lambda^{s,p}(\R^n)$ and $\Lambda^{s-t,p}(\R^n)$, for every $s,t\in\R$, $1\leq p\leq \infty$. In particular, $\Lambda^{s,p}(\R^n)$ is isomorphic to $L^p(\R^n)=\Lambda^{0,p}(\R^n)$. In fact, for any $u\in \Lambda^{s,p}(\R^n)$, by the semigroup property, $$\norm{\Lambda_{-t}u}_{\Lambda^{s-t,p}}=\norm{\Lambda_{t-s}(\Lambda_{-t}u)}_p=\norm{\Lambda_{t-s-t}u}_p=\norm{\Lambda_{-s}u}_p=\norm{u}_{\Lambda^{s,p}}.$$ 
This is usually known as the {\it lifting property of the Bessel potential}. We now define the complex interpolation spaces between $L^p$ spaces and $W^{k,p}$ spaces, which turns out to be equivalent to the spaces obtained by means of the Bessel potential.
\begin{Defi}\em[Complex interpolation of Sobolev spaces]
    Let $\theta\in [0,1]$, $1<p<\infty$ and $k\in \mathbb{N}$. Let $s=k\theta$. We define the space $H^{s,p}(\R^n)$ as the complex interpolation space $$H^{s,p}(\R^n):=[L^p(\R^n),W^{k,p}(\R^n)]_\theta,$$ with the identification $H^{0,p}(\R^n)=L^p(\R^n)$. In particular, for $s\in [0,1]$ and $k=1$,  
    $$H^{s,p}(\R^n)=[L^p(\R^n),W^{1,p}(\R^n)]_s.$$
\end{Defi}
\begin{Teor}\label{NormBessel}
    Let $\theta\in (0,1)$, $1<p<\infty$ and $k\in\mathbb{N}$. For $s=k\theta$, we have that $$\Lambda^{s,p}(\R^n)=H^{s,p}(\R^n),$$ with equivalence of norms. 
\end{Teor}
 We provide a direct proof of the equivalence of the interpolation norm and the norm in the Bessel potential space, which, to the authors' knowledge, is new in the literature. Part of the proof is based on some ideas of \cite[Chapter~13, Proposition~6.2]{Taylor2010}. In particular, the use of the following bound of the purely imaginary Bessel potentials, which follows directly from Mihlin's multiplier theorem \cite[Theorem~6.2.7]{Grafakos2014}.
\begin{Lema}
    Let $1<p<\infty,$ $u\in L^p(\R^n)$ and $t\in \R$. Then there exists a positive constant $C$ only depending on $p$ such that $$\norm{\Lambda_{it}u}_p\leq C\left(1+4\pi^2|t|^2\right)^{n/2}\norm{u}_p.$$
\end{Lema}
\noindent\textbf{Proof of Theorem \ref{NormBessel}:} We first prove the inclusion $$H^{s,p}(\R^n)\xhookleftarrow{}\Lambda^{s,p}(\R^n).$$ Let $u\in \Lambda^{s,p}(\R^n),$ we want to show that $u\in [L^p(\R^n),W^{k,p}(\R^n)]_\theta$, i.e., we must find a function $f\in \mathfrak{F}\left(L^p(\R^n),W^{k,p}(\R^n)\right)$ such that $f(\theta)=u$. Let $$f(z)=e^{(z-\theta)^2}\Lambda_{(z-\theta)k}u,$$ for $z\in S=\{z\in \C: 0\leq \operatorname{Re}z\leq 1\}$. Clearly, $f$ is a holomorphic function on $\operatorname{int}S$ and continuous in $S$, taking values on $L^p(\R^n)+W^{k,p}(\R^n)$. Now, we have to study the boundedness of $f(j+it)$ for $t\in \R$, $j=0,1$. We observe that $$f(j+it)=e^{-t^2+(\theta-j)^2+2i(jt-\theta t)}\Lambda_{(j+it-\theta)k}u.$$ Since we have that $$\left|e^{-t^2+(\theta-j)^2+2i(jt-\theta t)}\right|=e^{-t^2+(\theta-j)^2}\to 0,\,|t|\to \infty,$$ we only have to study the boundedness of 
$$\norm{\Lambda_{(j+it-\theta)k}u}_{W^{jk,p}},$$ where $W^{0,p}(\R^n)$ is identified with $L^p(\R^n)$. First, we have that 
$$\Lambda_{(it-\theta)k}u=\Lambda_{itk}\left(\Lambda_{-\theta k}u\right),$$ 
due to the semigroup properties of the Bessel potential. By the lifting property, $$\Lambda_{-\theta k}:\Lambda^{\theta k,p}(\R^n)\to \Lambda^{0,p}(\R^n)=L^p(\R^n),$$ hence $\Lambda_{-\theta k}u\in L^p(\R^n)$, and by the preceding lemma, 
\begin{align*}
    \norm{\Lambda_{(it-\theta)k}u}_p\leq C\left(1+4\pi^2|tk|^2\right)^{n/2}\norm{\Lambda_{-\theta k}u}_p<\infty,
\end{align*}
and hence 
$$\norm{f(it)}_p\leq Ce^{-t^2+\theta^2}\left(1+4\pi^2|tk|^2\right)^{n/2}\norm{\Lambda_{-\theta k}u}_p<\infty.$$ 
Observe that since $\left(1+4\pi^2|tk|^2\right)^{n/2}$ has at most polynomial growth in $t$, 
$$\lim_{|t|\to \infty}e^{-t^2+\theta^2}\left(1+4\pi^2|tk|^2\right)^{n/2}=0,$$ and hence $$\norm{f(it)}_p\to 0,$$ as $|t|\to \infty$. For $j=1$, we have that $$\Lambda_{(1+it-\theta)k}u=\Lambda_{itk}\left(\Lambda_k\left(\Lambda_{-\theta k}u\right)\right),$$
and since \[\begin{tikzcd}
	{\Lambda^{\theta k,p}(\R^n)} & {L^{p}(\R^n)} & {W^{k,p}(\R^n)}
	\arrow["{\Lambda_{-\theta k}}", from=1-1, to=1-2]
	\arrow["{\Lambda_k}", from=1-2, to=1-3]
\end{tikzcd}\] $\Lambda_k\left(\Lambda_{-\theta k}u\right)=\Lambda_{(1-\theta)k}u\in W^{k,p}(\R^n)$. Now, again by the preceding lemma and the fact that $W^{k,p}(\R^n)\xhookrightarrow{}L^p(\R^n)$ we get that $$\norm{\Lambda_{(1+it-\theta)k}u}_{W^{k,p}}\leq C'\left(1+4\pi^2|tk|^2\right)^{n/2}\norm{\Lambda_{(1-\theta)k}u}_{W^{k,p}}<\infty,$$ so $$\norm{f(1+it)}_{W^{k,p}}\leq C'e^{-t^2+(\theta-j)^2}\left(1+4\pi^2|tk|^2\right)^{n/2}\norm{\Lambda_{(1-\theta)k}u}_{W^{k,p}}<\infty,$$ which tends to zero again when $|t|\to \infty$. We conclude that $f\in \mathfrak{F}\left(L^p(\R^n),W^{k,p}(\R^n)\right)$, and since $$f(\theta)=\Lambda_0u=u,$$ we have that $u\in [L^p(\R^n),W^{k,p}(\R^n)]_\theta$, and hence $$\Lambda^{s,p}(\R^n)\xhookrightarrow{}H^{s,p}(\R^n).$$
Now we prove the reverse embedding. Let $v\in H^{s,p}(\R^n)$. Then, there exists $g\in \mathfrak{F}\left(L^p(\R^n),W^{k,p}(\R^n)\right),$ such that $g(\theta)=v$.  We have to prove that $g(\theta)\in \Lambda^{s,p}(\R^n)$. Since $g(\theta)\in L^p(\R^n)$, it is enough to prove that $\Lambda_{-s}g(\theta)\in L^p(\R^n)$. For every $\varphi\in \mathcal{S}(\R^n)$ and $z\in \overline{S}$, we define the bilinear form $$G_\varphi(z):=e^{(z-\theta)^2}\int_{\R^n}g(z)(x)\Lambda_{-kz}\varphi(x)\,dx,$$ where $g(z)(x)$ must be understood as the function $g(z)$ in $L^p(\R^n)$ for each $z\in \overline{S}$. Note that for each choice of $z$, this bilinear form defines a natural pairing for Schwartz functions, and for a fixed $\varphi$, it defines an holomorphic function in $S$ which is continuous and bounded up the boundary of the strip.

Let $M:=\norm{g}_{\mathfrak{F}\left(L^p,W^{k,p}\right)}$.
Now, let $z=it$ for $t\in \R$. By H\"older's inequality we get that $$|G_\varphi(it)|\leq e^{\theta^2-t^2}\norm{g(it)}_p\norm{\Lambda_{-kit}\varphi}_{q}\leq CMe^{\theta^2-t^2}(1+4\pi^2k^2t^2)^{n/2}\norm{\varphi}_q,$$ where $1/p+1/q=1$. Since $|G_\varphi(it)|\to 0$ as $|t|\to \infty$, we get that there exists a constant $C_0$ not depending on $t$ such that $$\operatorname{sup}_{t\in \R}|G_\varphi(it)|\leq C_0\norm{\varphi}_q.$$ Now, for $z=1+it$, we observe that by the semigroup property, $\Lambda_{-k(1+it)}=\Lambda_{-k}\Lambda_{-kit}$, and hence $$G_\varphi(1+it)=e^{(1-\theta)^2-t^2}\langle g(1+it),\Lambda_{-k(1+it)}\varphi\rangle=e^{(1-\theta)^2-t^2}\langle\Lambda_{-k}g(1+it),\Lambda_{-ikt}\varphi\rangle.$$
Observe that since $g(1+it)\in W^{k,p}(\R^n)=\Lambda^{k,p}(\R^n)$, we can use that $\Lambda_{-k}g(1+it)\in L^p(\R^n)$, and hence by an analogous reasoning we get that there exists a positive constant $C_1$ not depending on $t$ such that $$|G_{\varphi}(1+it)|\leq C_1\norm{\varphi}_q.$$
Now, by the Hadamaard's three-line theorem, $$|G_\varphi(\theta)|\leq C_0^{1-\theta}C_1^\theta\norm{\varphi}_{q},\,\forall \varphi\in \mathcal{S}(\R^n),$$ and by density we can extend this inequality for functions $\varphi\in L^q(\R^n)$. Finally, we observe that $$G_\varphi(\theta)=\langle g(\theta),\Lambda_{-s}\varphi\rangle=\langle \Lambda_{-s}v,\varphi\rangle,$$ so the functional $\varphi\mapsto G_\varphi(\theta)$, defines a bounded functional in $L^q(\R^n)$, and by duality it implies that we can see $\Lambda_{-s}v$ as a function on $L^p(\R^n)$ as we wanted to prove. \qed

We can extend the definition of complex interpolation spaces to negative exponents as $$H^{-k\theta,p}(\R^n):=[L^p(\R^n),W^{-k,p}(\R^n)]_\theta,\,k\in \mathbb{N}, p\in (1,\infty),\,\theta\in (0,1),$$ where $$W^{-k,p}(\R^n)=\left(W^{k,q}(\R^n)\right)^*,\,\frac{1}{p}+\frac{1}{q}=1.$$  By {Theorem \ref{ComplexProps}} we have the following natural result.
\begin{Prop}\label{BesselReflex}
    Let $s<0$, $1<p<\infty$ Then, $$\Lambda^{s,p}(\R^n)=H^{s,p}(\R^n),$$ with equivalence of the norms. Moreover, the space $H^{t,p}(\R^n)$ is reflexive for every $t\in\R$ and
\end{Prop}
\noindent\textbf{Proof:} Let $k\in\mathbb{N}$ and $\theta\in (0,1)$ be such that $s=k\theta$. Since $L^p(\R^n)\cap W^{k,p}(\R^n)=W^{k,p}(\R^n)$, the pair $\left(L^p(\R^n),W^{k,p}\right)$ is regular, hence by {Theorem \ref{ComplexProps} (6)} we have that \begin{align*}\left(H^{k\theta,p}(\R^n)\right)^*&=\left([L^p(\R^n),W^{k,p}(\R^n)]\right)^*=\left[\left(L^p(\R^n)\right)^*,\left(W^{k,p}(\R^n)\right)^*\right]_\theta\\
&=[L^q(\R^n),W^{-k,q}(\R^n)]_\theta=H^{-k\theta,q}(\R^n),\,\frac{1}{p}+\frac{1}{q}=1.\end{align*} 
Now, since by {Theorem \ref{NormBessel}} $$\Lambda^{k\theta,p}(\R^n)=H^{k\theta,p}(\R^n),$$ 
it only remains to prove that $\Lambda^{-k\theta,q}(\R^n)=\left(\Lambda^{k\theta,p}(\R^n)\right)^*$, since $\left(\Lambda^{k\theta,p}(\R^n)\right)^*=\left(H^{s,p}(\R^n)\right)^*=H^{-k\theta,q}(\R^n)$.

Let $g\in \Lambda^{-k\theta,q}(\R^n)$, hence there exists $g_{-k\theta}\in L^q(\R^n)$ such that $\Lambda_{k\theta}g=g_{-k\theta}$, so $g=\Lambda_{-k\theta}g_{-k\theta}$. Now take $f\in \mathcal{S}(\R^n)$ and $f_{k\theta}=\Lambda_{-k\theta}f\in \Lambda^{k\theta,p}(\R^n)$. Now, in the sense of pairing of distributions and continuous linear functionals, 
\begin{align*}
    \langle g,f\rangle&=\langle \Lambda_{-k\theta}g_{-k\theta},\Lambda_{k\theta}f_{k\theta}\rangle=\left\langle \mathfrak{F}^{-1}\left((1+4\pi^2|\xi^2|)^{k\theta/2}\mathfrak{F}g_{-k\theta}(\xi)\right),\mathfrak{F}^{-1}\left((1+4\pi^2|\xi^2|)^{-k\theta/2}\mathfrak{F}f_{k\theta}(\xi)\right)\right\rangle\\
    &=\left\langle (1+4\pi^2|\xi^2|)^{k\theta/2}\mathfrak{F}g_{-k\theta}(\xi),(1+4\pi^2|\xi^2|)^{-k\theta/2}\mathfrak{F}f_{k\theta}(\xi)\right\rangle\\
    &=\left\langle\mathfrak{F}g_{-k\theta}(\xi),\mathfrak{F}f_{k\theta}(\xi)\right\rangle
    =\langle g_{-k\theta}(\xi),f_{k\theta}(\xi)\rangle\leq \norm{g_{-k\theta}}_q\norm{f_{k\theta}}_p\\
    &=\norm{g_{-k\theta}}_q\norm{\Lambda_{-k\theta}f}_p=\norm{g_{-k\theta}}_q\norm{f}_{\Lambda^{k\theta,p}}.
\end{align*}
Now, by {Proposition \ref{SchwartzDenseBessel}}, we can extend this result for functions $f\in \Lambda^{k\theta,p}(\R^n)$, and hence we can identify $g$ with a continuous linear functional on $\Lambda^{k\theta,p}(\R^n)$, which means that $$\Lambda^{-k\theta,q}(\R^n)\xhookrightarrow{}\left(\Lambda^{k\theta,p}(\R^n)\right)^*.$$
Now, consider $v\in \left(\Lambda^{k\theta,p}(\R^n)\right)^*\xhookrightarrow{}\mathcal{S}'(\R^n)$, so we can regard $v$ as a tempered distribution and hence its Fourier transform is well defined. Let $u\in \mathcal{S}(\R^n)$, hence $\Lambda_{-k\theta}u\in \mathcal{S}(\R^n)$, and we define $u_{k\theta}=\Lambda_{-k\theta}u\in \mathcal{S}(\R^n)$. Again, in the sense of pairings of distributions, \begin{align*}
    \langle \Lambda_{k\theta}v,u_{k\theta}\rangle&=\left\langle\mathfrak{F}^{-1}\left((1+4\pi^2|\xi|^2)^{-k\theta/2}\mathfrak{F}v\right),u_{k\theta}\right\rangle=\left\langle (1+4\pi^2|\xi|^2)^{-k\theta/2}\mathfrak{F}v,\mathfrak{F}^{-1}u_{k\theta}(\xi)\right\rangle\\
    &=\left\langle\mathfrak{F}v,(1+4\pi^2|\xi|^2)^{-k\theta/2},\mathfrak{F}u_{k\theta}(-\xi)\right\rangle=\left\langle v,\mathfrak{F}\left((1+4\pi^2|\xi|^2)^{-k\theta/2}\mathfrak{F}u_{k\theta}(-\xi)\right)\right\rangle\\
    &=\left\langle v,\mathfrak{F}^{-1}\left((1+4\pi^2|\xi|^2)^{-k\theta/2}\mathfrak{F}u_{k\theta}(\xi)\right)\right\rangle=\langle v,\Lambda_{k\theta}u_{k\theta}\rangle\leq\norm{v}_{\left(\Lambda^{k\theta,p}(\R^n)\right)^*}\norm{\Lambda_{k\theta}u_{k\theta}}_{\Lambda^{k\theta,p}}\\
    &= \norm{v}_{\left(\Lambda^{k\theta,p}(\R^n)\right)^*}\norm{u_{k\theta}}_p.
\end{align*}
Since $\mathcal{S}(\R^n)$ is dense in $L^p(\R^n)$ we can extend the result to functions in $L^p(\R^n)$, and hence $\Lambda_{k\theta}v\in \left(L^p(\R^n)\right)^*=L^q(\R^n)$, i.e., $v\in \Lambda^{-k\theta,q}(\R^n)$, so $$\left(\Lambda^{k\theta,p}(\R^n)\right)^*\xhookrightarrow{}\Lambda^{-k\theta,q}(\R^n).$$
Finally, take $m\in \mathbb{Z}$. Since both $L^p(\R^n)$ and $W^{m,p}(\R^n)$ are reflexive, by {Theorem \ref{ComplexProps} (7)}, the space $H^{m\theta,p}(\R^n)$ is reflexive.\qed\\

Henceforth, we shall denote $H^{s,p}(\mathbb{R}^n)$ in place of $\Lambda^{s,p}(\mathbb{R}^n)$ and refer to these as Bessel potential spaces, or simply Bessel spaces, for all $s \in \mathbb{R}$ and $1 < p < \infty$. We now summarize some fundamental properties of these spaces, which can be directly derived from established abstract interpolation results. 

\begin{Prop}
    For every $s\in\R$ and $1<p<\infty$, the space $H^{s,p}(\R^n)$ is complete.
\end{Prop}
\noindent\textbf{Proof:} For every compatible couple of Banach spaces $(E_0,E_1)$, the complex interpolation space $[E_0,E_1]_{\theta}$, $\theta\in (0,1)$ is complete since $\Cc_\theta$ is an interpolation functor, hence the result follows from the completeness of the spaces $L^p(\R^n), W^{k,p}(\R^n)$, $k\in\mathbb{N}$. \qed
\begin{Prop}\label{SobodenseBessel}
    Let $1<p<\infty$, $k\in \mathbb{Z}$ and $\theta\in (0,1)$. Then $W^{k,p}(\R^n)$ is dense $H^{k\theta,p}(\R^n)$.
\end{Prop}
\noindent\textbf{Proof:} Given any compatible couple of Banach spaces $(E_0,E_1)$, the space $E_0\cap E_1$ is dense in $[E_0,E_1]_\theta$ by {Proposition \ref{ComplexProps} (4)}. Since $L^p(\R^n)\cap W^{k,p}(\R^n)=W^{k,p}(\R^n)$, it follows that $W^{k,p}(\R^n)$ is dense in $H^{k\theta,p}(\R^n)$.\qed

\begin{Prop}\label{BesselAnidados}
Let $0<t<s$ and $1<p<\infty$. Then, we have the inclusion $$H^{s,p}(\R^n)\xhookrightarrow{}H^{t,p}(\R^n).$$
\end{Prop}
\noindent\textbf{Proof:}
    Let $k$ the least integer greater than $s$. By {Proposition \ref{ComplexProps} (2)} we know that if $E_0\xhookrightarrow{}E_1$ and $0<\theta_0<\theta_1<1$, then $$[E_0,E_1]_{\theta_0}\xhookrightarrow{}[E_0,E_1]_{\theta_1}.$$ Now, since $W^{k,p}(\R^n)\xhookrightarrow{}L^p(\R^n)$, and $1-s/k<1-t/k$, we get that \begin{align*}
        H^{s,p}(\R^n)&=[L^p(\R^n),W^{k,p}(\R^n)]_{s/k}=[W^{k,p}(\R^n),L^p(\R^n)]_{1-s/k}\xhookrightarrow{}[W^{k,p}(\R^n),L^p(\R^n)]_{1-t/k}\\
        &=[L^p(\R^n),W^{k,p}(\R^n)]_{t/k}=H^{t,p}(\R^n),
    \end{align*}
    where the equality $$[L^p(\R^n),W^{1,p}(\R^n)]_{s/k}=[W^{k,p}(\R^n),L^p(\R^n)]_{1-s/k},$$ follows from {Proposition \ref{ComplexProps} (1)}. \qed
    \begin{Prop}
        Let $s=k\in\mathbb{N}$. Then, the spaces $H^{s,p}(\R^n)$ coincide with the classical Sobolev spaces $W^{k,p}(\R^n)$ for every $1<p<\infty$.
    \end{Prop}
    \noindent\textbf{Proof:} For any compatible couple of Banach spaces $(E_0,E_1)$, by {Proposition \ref{ComplexProps} (5)} we have that $[E_0,E_1]_j$ is a closed subspace of $E_j$ with coincidence of the norms in $[E_0,E_1]_j$, $j=0,1$. In this case, it implies that $$H^{k,p}(\R^n)=[L^p(\R^n),W^{k,p}(\R^n)]_{1}\xhookrightarrow{}W^{k,p}(\R^n),$$ and since $E_0\cap E_1\xhookrightarrow{}[E_0,E_1]_\theta$, for all $\theta\in[0,1]$,  
    $$W^{k,p}(\R^n)=L^p(\R^n)\cap W^{k,p}(\R^n)\xhookrightarrow{}[L^p(\R^n),W^{k,p}(\R^n)]_1=H^{k,p}(\R^n),$$ 
    so $H^{k,p}(\R^n)=W^{k,p}(\R^n)$ with equivalence of the norms.\qed
    \begin{Prop}
        Let $s>0$ and $1<p<\infty$. Then, the space $C_c^\infty(\R^n)$ is dense in $H^{s,p}(\R^n)$.
    \end{Prop}
    \noindent\textbf{Proof:} Suppose that $s=k\theta$ for some $k\in\mathbb{N}$ and $\theta\in (0,1)$. Let $\varepsilon>0$ and $u\in H^{k\theta,p}(\R^n)$. By {Proposition \ref{SobodenseBessel}}, the space $W^{k,p}(\R^n)$ is dense in $H^{k\theta,p}(\R^n)$, and additionally 
    $$\norm{\cdot}_{H^{s,p}}\leq \norm{\cdot}_{W^{k,p}}.$$ 
    Hence, there exists $v\in W^{k,p}(\R^n)$ such that $$\norm{u-v}_{H^{s,p}}<\varepsilon/2.$$ 
    Also, since $C_c^\infty(\R^n)$ is dense in $W^{k,p}(\R^n)$, there exists $w\in C_c^\infty(\R^n)$ such that $$\norm{v-w}_{W^{k,p}(\R^n)}<\varepsilon/2.$$ Then, we have $$\norm{u-w}_{H^{s,p}}\leq \norm{u-v}_{H^{s,p}}+\norm{v-w}_{H^{s,p}}\leq \norm{u-w}_{H^{s,p}}+\norm{v-w}_{W^{k,p}}<\varepsilon/2+\varepsilon/2=\varepsilon,$$ so $C_c^\infty(\R^n)$ is dense in $H^{s,p}(\R^n)$.\qed
    \begin{Prop}
        Let $s\in \R$ and $1<p<\infty$. Then, the space $H^{s,p}(\R^n)$ is separable.
    \end{Prop}
    \noindent\textbf{Proof:} For $s>0$, the result follows from the fact that subspaces of separable metric spaces are separable. Indeed, $L^p(\R^n)$ is separable for $1\leq p<\infty$ and $H^{s,p}(\R^n)\xhookrightarrow{}L^p(\R^n)+W^{1,p}(\R^n)=L^p(\R^n)$, we get the desired result. Now, for $s<0$, the result follows from the fact that Banach spaces with separable dual are separable. Since by {Proposition \ref{BesselReflex}} $$H^{s,p}(\R^n)=\left(H^{-s,q}(\R^n)\right)^*,\,\frac{1}{p}+\frac{1}{q}=1,$$ and Bessel potential spaces are reflexive, and since $H^{-s,q}(\R^n)=\left(H^{-s,q}(\R^n)\right)^{**}=\left(H^{s,p}(\R^n)\right)^*$ is separable, then $H^{s,p}(\R^n)$ for $s<0$ is separable. For $s=0$, separability follows from the identification $H^{0,p}(\R^n)=L^p(\R^n)$. \qed\\
    
    Now, we are going to prove several further results using finer results on the complex and real interpolation methods.
    \begin{Teor}[Real reiteration of Bessel spaces] Let $s\in \R$, $1<p_0,p_1<\infty$. Then, $$\left(H^{s,p_0}(\R^n),H^{s,p_1}(\R^n)\right)_{\theta,p(\theta)}=H^{s,p(\theta)}(\R^n),\,\frac{1}{p(\theta)}=\frac{1-\theta}{p_0}+\frac{\theta}{p_1}.$$
    \end{Teor}
    \noindent\textbf{Proof:} Since $\Lambda_{-t}$ is an isomorphism between $H^{s,p}(\R^n)$ and $H^{s-t,p}(\R^n)$, $\Lambda_{-s}$ is an isomorphism between $H^{s,p}(\R^n)$ and $H^{0,p}(\R^n)$, which is identified with $L^p(\R^n)$. Then, using this and the results on real interpolation of Lebesgue spaces,
    \begin{align*}
        \left(H^{s,p_0}(\R^n),H^{s,p_1}(\R^n)\right)_{\theta,p(\theta)}&=\left(\Lambda_{-s}H^{s,p_0}(\R^n),\Lambda_{-s}H^{s,p_1}(\R^n)\right)_{\theta,p(\theta)}=\left(L^{p_0}(\R^n),L^{p_1}(\R^n)\right)_{\theta,p(\theta)}\\
        &=L^{p(\theta),p(\theta)}(\R^n)=L^{p(\theta)}(\R^n)=\Lambda_{-s}H^{s,p(\theta)}(\R^n)=H^{s,p}(\R^n),
    \end{align*}
    where the equalities are up to the equivalence of norms. \qed\\
  
For the reiteration using the complex method, it is necessary to employ the concept of retraction and {Theorem \ref{retract}}. The key is that the spaces $H^{s,p}(\mathbb{R}^n)$ are retracts of the Lebesgue spaces of the weighted sequences $L^p\left(\ell^{s,2}\right)$ for $s \in \mathbb{R}$ and $1 < p < \infty$. This result is a direct consequence of Mihlin's multiplier theorem. Since it involves the Littlewood-Paley characterization of Bessel spaces, which is beyond the scope of this paper, we redirect the interested reader to \cite[Theorem~6.4.3]{BerghLofstrom1976}, which we rely on.

    \begin{Teor}[Complex reiteration of Bessel spaces]\label{complexreiter}
       Let $0<\theta<1$, $s_0,s_1\in \R$, $s_0\not =s_1$, $1<p_0,p_1<\infty$. Then, 
       $$[H^{s_0,p_0}(\R^n),H^{s_1,p_1}(\R^n)]_{\theta}=H^{s(\theta),p(\theta)},\,s(\theta)=(1-\theta)s_0+\theta s_1,\,\frac{1}{p(\theta)}=\frac{1-\theta}{p_0}+\frac{\theta}{p_1}.$$
    \end{Teor}
    \noindent\textbf{Proof:} The case $p_0=p_1$ follows easily from {Theorem \ref{ComplexReiter}}. If $s_0,s_1\in (0,1)$, \begin{align*}[H^{s_0,p}(\R^n),H^{s_1,p}(\R^n)]_\theta&=\left[[L^p(\R^n),W^{1,p}(\R^n)]_{s_0},[L^p(\R^n),W^{1,p}(\R^n)]_{s_1}\right]_\theta\\
    &=[L^p(\R^n),W^{1,p}(\R^n)]_{s(\theta)}=H^{s(\theta),p}(\R^n).\end{align*} If $|s_0|,|s_1|>1$, we obtain the result just interpolating with respect to the couple $\left(L^p(\R^n),W^{k,p}(\R^n)\right)$, with $k$ being the least integer bigger than $|s_0|$ and $|s_1|$, and the interpolation parameters $|s_i'|=|s_i|/k$, $i=0,1$.\\
    Now, for the general case $p_0\not =p_1$ we have to use the fact that $H^{s,p}(\R^n)$ is a retract of $L^p\left(\ell^{s,2}\right)$. Indeed, {Theorem \ref{retract}} and the fact that for any compatible couple of Banach spaces $(E_0,E_1)$, 
    $$[L^{p_0}(E_0),L^{p_1}(E_1)]_\theta=L^{p(\theta)}\left([E_0,E_1]_\theta\right),$$ 
    for any $1\leq p_0,p_1<\infty$ and $\theta\in (0,1)$. Also, since  $$[\ell^{s_0,q_0},\ell^{s_1,q_1}]_\theta=\ell^{s(\theta),q(\theta)},$$ for any $s_0,s_1\in \R$ and $1\leq q_0,q_1\leq \infty$, we have that $$[L^{p_0}\left(\ell^{s_0,2}\right),L^{p_1}\left(\ell^{s_1,2}\right)]_\theta=L^{p(\theta)}\left(\left[\ell^{s_0,2},\ell^{s_1,2}\right]_\theta\right)=L^{p(\theta)}\left(\ell^{s(\theta),2}\right),$$ with the latest space being a coretract of $H^{s(\theta),p(\theta)}(\R^n)$, so we conclude that 
    \[
\pushQED{\qed} 
[H^{s_0,p_1}(\R^n),H^{s_1,p_1}(\R^n)]_\theta=H^{s(\theta),p(\theta)}(\R^n).\qedhere
\popQED
\] 
    \begin{Obse}\em
        The last theorem does not cover the cases $p=1$ and $p=\infty$, since they are based on the reiteration theorem for the complex method and Mihlin's multiplier theorem, as $H^{s,1}(\R^n)$ and $H^{s,\infty}(\R^n)$ are not covered by those results. However, if we fix $s\in \R$, by the lifting property of the Bessel potential, we have that $$[H^{s,1}(\R^n),H^{s,\infty}(\R^n)]_{\theta}=\Lambda_s [L^1(\R^n),L^\infty(\R^n)]_\theta=\Lambda_s L^{\frac{1}{1-\theta}}(\R^n)=H^{s,\frac{1}{1-\theta}}(\R^n).$$ Now, by {Theorem \ref{complexreiter}}, given $1<p<\infty$, \begin{align*}[H^{s,1}(\R^n),H^{s,p}(\R^n)]_\theta&=\left[H^{s,1}(\R^n),[H^{s,1}(\R^n),H^{s,\infty}(\R^n)]_{\frac{p-1}{p}}\right]_\theta=[H^{s,1}(\R^n),H^{s,\infty}(\R^n)]_{\theta\frac{p-1}{p}}\\
        &=H^{s,p(\theta)}(\R^n),\end{align*} where $$\frac{1}{p(\theta)}=1-\theta\frac{p-1}{p}=(1-\theta)+\frac{\theta}{p}.$$ 
        And in a similar fashion, 
        \begin{align*}
            [H^{s,p}(\R^n),H^{s,\infty}(\R^n)]_\theta&=\left[[H^{s,1}(\R^n),H^{s,\infty}(\R^n)]_{\frac{p-1}{p}},H^{s,\infty}(\R^n)\right]_\theta=[H^{s,1}(\R^n),H^{s,\infty}(\R^n)]_{(1-\theta)\frac{p-1}{p}+\theta}\\
        &=H^{s,\tilde{p}(\theta)}(\R^n),
        \end{align*}
        where $$\frac{1}{\tilde{p}(\theta)}=1-(1-\theta)\frac{p-1}{p}-\theta=\frac{1-\theta}{p},$$
        extending the previous result to the cases $p=1$ and $p=\infty$. 
    \end{Obse}
    This should be compared with the results of Milman \cite{Milman1983}, where it was established that for $k\in \mathbb{N}$ and $1<p<\infty$, $$[W^{k,1}(\R^n),W^{k,p}(\R^n)]_\theta=W^{k,p(\theta)}(\R^n),\,\frac{1}{p(\theta)}=1-\theta+\frac{\theta}{p}.$$ 
    However, there is a notable difference with the fractional case since the techniques of Milman do not extend to the case $p=\infty$. Recently, it was finally proved by Curca \cite{curca2024} that there is no possible extension to $p=\infty$ in the case of classical Sobolev spaces, solving the long-standing question posed in 1984 by Jones (see \cite{jones1984}). This differs from the real method case, for which it was soon established by De Vore and Scherer \cite{devore1979interpolation} that it covers the case $p=\infty$ 
\begin{Coro}[Interpolation inequalities for Bessel potential spaces]
    Let $s_0\not =s_1\in \R$, $1<p_0,p_1<\infty$ and $\theta\in (0,1)$. Then, there exists a positive constant $C$ such that for every $ u \in H^{s_0,p_0}\cap H^{s_1,p_1}$ $$\norm{u}_{H^{s(\theta),p(\theta)}}\leq C\norm{u}_{H^{s_0,p_0}}^{1-\theta}\norm{u}_{H^{s_1,p_1}}^\theta.$$
\end{Coro}
\noindent\textbf{Proof:} From the complex reiteration theorem $$[H^{s_0,p_0}(\R^n),H^{s_1,p_1}(\R^n)]_{\theta}=H^{s(\theta),p(\theta)}(\R^n),$$ and since the functor $\Cc_\theta$ is an interpolation functor of exponent $\theta$, the result follows. \qed\\

Notice that in the previous result the constant $C$ can be chosen so that it depends continuously on the parameters.

From {Remark \ref{CoincidenciaDeMetodos}} it also follows that for every $s\in (0,1)$ and $1<p<\infty$ $$\left(L^p(\R^n),W^{1,p}(\R^n)\right)_{s,1}\xhookrightarrow{}H^{s,p}(\R^n)\xhookrightarrow{}\left(L^p(\R^n),W^{1,p}(\R^n)\right)_{s,\infty}.$$ The interpolation spaces at the endpoints are identified with the large class of Besov spaces \cite[Theorem~17.24]{Leoni2017}.

\subsection{Continuous Embeddings via Interpolation Theory}
In this section, we will establish the Sobolev embeddings in the fractional case using the interpolation of linear operators and the well-known results for classical Sobolev spaces. Some of these embeddings have been known since the beginning of the modern theory of function spaces and were initially established for the general three-parameter family of \textit{Triebel-Lizorkin} spaces \( F_{p,q}^s \), where \( 1 \leq p, q \leq \infty \) and \( s \in \mathbb{R} \). These spaces are connected to \( H^{s,p}(\mathbb{R}^n) \) via the Littlewood-Paley theory. Specifically, \( F^{s}_{p,2}(\mathbb{R}^n) = H^{s,p}(\mathbb{R}^n) \) for \( s \in \mathbb{R} \) and \( 1 < p < \infty \). However, the proofs (see \cite[Theorem~2.8.1]{Triebel1995}, for instance) often rely on indirect arguments, such as the Littlewood-Paley characterization of Bessel potential spaces. More recently, in \cite{ShiehSpector2015}, these embeddings were established using the Riesz fractional gradient and well-known estimates for the Riesz potential. Our proofs are entirely based on the classical embeddings for the endpoints of the interpolation scale and the properties of the complex method, making them particularly interesting from a qualitative perspective.\\

The first one is the embedding of $H^{s,p}(\R^n)$ into the space $\textbf{BMO}(\R^n)$ when $sp=n$. This result is of interest in its own right, and will be longely used in the proof of the Fractional Sobolev embedding.
\begin{Teor}\label{CriticalI} Let $s\in (0,1)$ and $1<p<\infty$ such that $sp=n$. Then, we have $$H^{s,p}(\R^n)\xhookrightarrow{}{\normalfont \textbf{BMO}}(\R^n).$$
\end{Teor}
\noindent\textbf{Proof:} Consider the operator $$T:L^p(\R^n)+W^{1,p}(\R^n)\to L^{p}(\R^n)+C^{0,\mu}(\R^n); u\mapsto u,$$ where $\mu=1-n/p$. This mapping is well defined and is an admissible operator, since $T|_{L^p(\R^n)}$ is just the identity map for $L^p(\R^n)$ and $T|_{W^{1,p}(\R^n)}\in \mathcal{L}\left(W^{1,p}(\R^n),C^{0,\mu}(\R^n)\right)$. As $p>n$ (because $sp=n$), Morrey's inequality holds (see \cite[Theorem~9.12]{Brezis2010}), thus we have that $W^{1,p}(\R^n)\xhookrightarrow{}C^{0,1-n/p}(\R^n)$. Hence, $T: u\mapsto u$ is a bounded linear map (in fact, it is an inclusion operator) from $[L^p(\R^n),W^{1,p}(\R^n)]_s$ to $[L^p(\R^n), C^{0,\mu}(\R^n)]_s$. The couple $\left(L^p(\R^n), C^{0,\mu}(\R^n)\right)$ can be identified with a compatible couple of \textit{Morrey-Campanato} spaces $\mathcal{L}^{p,\lambda}(\R^n)$, $1\leq p<\infty$, $0\leq \lambda\leq n+p$, (see \cite{Peetre1969II})
the space of locally $L^p(\R^n)$ functions $f$ such that $$|f|_{\mathcal{L}^{p,\lambda}(\R^n)}:=\operatorname{sup}_{x_0\in\R^n, r>0}r^{-\lambda/p}\left(\int_{B(x_0,r)}|f(x)-f_{B(x_0,r)}|^p\,dx\right)^{1/p}<\infty.$$
Clearly, $\mathcal{L}^{p,0}(\R^n)$ is equivalent to the space $L^p(\R^n)$, for $\lambda=n$ it is equivalent to $\textbf{BMO}(\R^n)$ and for $\lambda>n$, $\mathcal{L}^{p,\lambda}(\R^n)$ is equivalent to $C^{0,\alpha}(\R^n)$, where $$\alpha=\frac{\lambda-n}{p},$$ (see \cite[Theorem~7.1]{Leoni2023}). $C^{0,\mu}(\R^n)$ is equivalent to $\mathcal{L}^{p,\lambda}$ with $$\lambda=p\mu+n=p\left(1-n/p\right)+n=p,$$ and hence the couple $\left(L^p(\R^n), C^{0,\mu}(\R^n)\right)$ is equivalent to the couple $\left(\mathcal{L}^{p,0}(\R^n),\mathcal{L}^{p,p}(\R^n)\right)$, for which interpolation properties are well known (see the introduction in \cite{BlascoRuizVega1999}). Then, from the previous discussion $$H^{s,p}(\R^n)=[L^p(\R^n),W^{1,p}(\R^n)]_s\xhookrightarrow{}[\mathcal{L}^{p,0}(\R^n),\mathcal{L}^{p,p}(\R^n)]_s\xhookrightarrow{}\mathcal{L}^{p(s),\lambda(s)}(\R^n),$$ 
where 
$$\frac{1}{p(s)}=\frac{1-s}{p}+\frac{s}{p}=\frac{1}{p},$$
so
$$p(s)=p,$$ and $$\lambda(s)=(1-s)0+sp=sp.$$ Since $sp=n$, we have that $\mathcal{L}^{p,sp}(\R^n)$ is equivalent to the space $\textbf{BMO}(\R^n)$, so we conclude that \[\pushQED{\qed} 
H^{s,p}(\R^n)\xhookrightarrow{}\textbf{BMO}(\R^n).\qedhere
\popQED\]
\begin{Teor}[Fractional Sobolev embedding theorem]\label{FSET}
    Let $s\in (0,1)$ and $1<p<\infty$. Then: \begin{enumerate}
        \item\textbf{Subcritical case.-} If $sp<n$, then $$H^{s,p}(\R^n)\xhookrightarrow{}L^q(\R^n),\,q\in [p,p_s^*],$$ where $$p_s^*=\frac{np}{n-sp},$$ is the fractional Sobolev conjugate exponent.
        \item\textbf{Critical case.-} If $sp=n$, then $$H^{s,p}(\R^n)\xhookrightarrow{}L^q(\R^n),\,q\in [p,\infty).$$
        \item\textbf{Supercritical case or Morrey's embedding.-} If $sp>n$, then $$H^{s,p}(\R^n)\xhookrightarrow{}C^{0,\mu}(\R^n),\,0<\mu\leq \mu_s^*,$$
        with $\mu_s^*=s-\frac{n}{p}$.
    \end{enumerate}
\end{Teor}
\noindent\textbf{Proof:}
\begin{enumerate}
    \item We will study three different cases according to the relationship of $n$ and $p$. 
\begin{itemize}
    \item If $n>p$, 
    consider the operator $$T:L^p(\R^n)+W^{1,p}(\R^n)\xhookrightarrow{}L^p(\R^n)+L^{p^*}(\R^n); u\mapsto u,$$ where $$p^*=\frac{np}{n-p},$$ is well defined. Clearly $T\in \mathcal{L}\left(L^p(\R^n),L^p(\R^n)\right)$ since $T_{|_{L^p(\R^n)}}$ is just the identity map. Also,
$T\in \mathcal{L}\left(W^{1,p}(\R^n),L^{p^*}(\R^n)\right)$ since  it is well known that $W^{1,p}(\R^n)\xhookrightarrow{} L^{p^*}(\R^n)$ (see \cite[Theorem~9.9]{Brezis2010}, with norm less or equal to some constant $C$ depending only on $p$ and $n$. Hence by $T: u\mapsto u$ is an admissible operator from $\left(L^p(\R^n),W^{1,p}(\R^n)\right)$, and since $\Cc_\theta$ is an exact interpolation functor, $T$ is a bounded linear map (in fact it is an inclusion operator) from $[L^p(\R^n),W^{1,p}(\R^n)]_s$ to $[L^p(\R^n),L^{p^*}(\R^n)]_s$, and since 
$$[L^p(\R^n),W^{1,p}(\R^n)]_s=H^{s,p}(\R^n),\,[L^p(\R^n),L^{p^*}(\R^n)]_s=L^{p(s)}(\R^n),$$ 
where 
$$\frac{1}{p(s)}=\frac{1-s}{p}+\frac{s}{p^*}$$
implies
$$p(s)=\frac{np}{n-sp}=p_s^*,$$ we have that $$\norm{u}_{p_s^*}\leq C'\norm{u}_{H^{s,p}},$$ for some constant $C'$ depending on $n,p$ and $s$, thus 
$$H^{s,p}(\R^n)\xhookrightarrow{}L^{p_s^*}(\R^n).$$
\item If $p=n$, proceeding as in the previous case, interpolating the operator $$T:L^p(\R^n)+W^{1,p}(\R^n)\xhookrightarrow{} L^p(\R^n)+L^{p^*}(\R^n); u\mapsto u,$$ we would get that $H^{s,p}(\R^n)\xhookrightarrow{}L^{p_s^*}(\R^n)$. However, note that $p^*=\infty$, and it is known that the embedding $W^{1,n}(\R^n)\xhookrightarrow{}L^\infty (\R^n)$ fails to hold, so we need a suitable space replacing $L^\infty(\R^n)$. Consider the space $\textbf{BMO}(\R^n)$ ({Examples \ref{ComplexEjemes}} iii)). We know that $$[L^p(\R^n),\textbf{BMO}(\R^n)]_s=L^{q(s)}(\R^n),$$ where $\frac{1}{q(s)}=\frac{1-s}{p}$, and since $p=n$, $$p_s^*=\frac{np}{n-sp}=\frac{np}{n-sn}=\frac{p}{1-s}=q(s),$$ so $$[L^p(\R^n),\textbf{BMO}(\R^n)]_s=L^{p_s^*}(\R^n).$$ Now, since $W^{1,n}(\R^n)\xhookrightarrow{}\textbf{BMO}(\R^n)$ (see \cite{BrezisNirenberg1995}), the operator $$T:L^p(\R^n)+W^{1,n}(\R^n)\xhookrightarrow{}L^p(\R^n)+\textbf{BMO}(\R^n),$$ is an admissible operator (in fact, it is an inclusion operator), and therefore we conclude that $$H^{s,p}(\R^n)=[L^p(\R^n),W^{1,p}(\R^n)]_s\xhookrightarrow{}[L^p(\R^n),\textbf{BMO}(\R^n)]_s=L^{p_s^*}(\R^n).$$
\item If $p>n$, there exists $t\in (s,1)$ such that $tp=n$. Now, consider the operator $$T:L^p(\R^n)+H^{t,p}(\R^n)\to L^p(\R^n)+\textbf{BMO}(\R^n),$$ which is an admissible inclusion operator, since by {Theorem \ref{CriticalI}} $$H^{t,p}(\R^n)\xhookrightarrow{}\textbf{BMO}(\R^n).$$ Then, $T$ is a bounded linear operator (inclusion) between $[L^p(\R^n),H^{t,p}(\R^n)]_{s/t}$ and\\
$[L^p(\R^n),\textbf{BMO}(\R^n)]_{s/t}$. By {Theorem \ref{complexreiter}}, $$[L^p(\R^n),H^{t,p}(\R^n)]_{s/t}=H^{s,p}(\R^n),$$ and since, $$[L^p(\R^n),\textbf{BMO}(\R^n)]_{s/t}=L^q(\R^n),$$ where 
$$\frac{1}{q}=\frac{1-s/t}{p} $$ 
implies
$$q=\frac{p}{1-\frac{s}{t}}=\frac{np}{n-n\frac{s}{t}}.$$ Note that $tp=n$, hence $$q=\frac{np}{n-tp\frac{s}{t}}=\frac{np}{n-sp}=p_s^*,$$ thus $$H^{s,p}(\R^n)\xhookrightarrow{}L^{p_s^*}(\R^n).$$ 
\end{itemize}
Now, let $p<q<p_s^*$. Since $H^{s,p}(\R^n)\xhookrightarrow{}L^p(\R^n)$, we have that $$H^{s,p}(\R^n)\xhookrightarrow{}[L^p(\R^n),L^{p_s^*}(\R^n)]_\theta,$$ for any $0<\theta<1$. Taking $$\theta=\frac{p_s^*(q-p)}{q(p_s^*-p)}$$ yields $$H^{s,p}(\R^n)\xhookrightarrow{}L^q(\R^n).$$
\item By {Proposition \ref{BesselAnidados}}, $$H^{s,p}(\R^n)\xhookrightarrow{}H^{s-\varepsilon,p}(\R^n),$$ for every $\varepsilon\in (0,s)$. Since $sp=n$, we have that $$0<(s-\varepsilon)p=sp-\varepsilon p=n-\varepsilon p<n,$$ hence by the subcritical case, $$H^{s-\varepsilon,p}(\R^n)\xhookrightarrow{}L^{p_{s-\varepsilon}^*}(\R^n),$$ where $$p_{s-\varepsilon}^*=\frac{np}{n-(s-\varepsilon)p}.$$ Since the function $\varepsilon\mapsto p_{s-\varepsilon}^*$ is a continuous function mapping $(0,s)$ to $(p,\infty)$, we have that $$H^{s,p}(\R^n)\xhookrightarrow{}L^q(\R^n),$$ for every $q\in (p,\infty)$. The case $p=q$ is trivial since $H^{s,p}(\R^n)$ is a subspace of $L^p(\R^n)$.
\item Analogously as in the proof of {Theorem \ref{CriticalI}}, we get that $$H^{s,p}(\R^n)=[L^p(\R^n),W^{1,p}(\R^n)]_s\xhookrightarrow{}[\mathcal{L}^{p,0}(\R^n),\mathcal{L}^{p,p}(\R^n)]_s\xhookrightarrow{}\mathcal{L}^{p(s),\lambda(s)}(\R^n),$$ 
where 
$$\frac{1}{p(s)}=\frac{1-s}{p}+\frac{s}{p}=\frac{1}{p}$$
implies
$$p(s)=p,$$ and $$\lambda(s)=(1-s)0+sp=sp.$$ 
Since $sp>n$, we have that $\mathcal{L}^{p,sp}(\R^n)$ is equivalent to the space $C^{0,\mu_s^*}(\R^n)$, so we conclude that $$H^{s,p}(\R^n)\xhookrightarrow{}C^{0,s-n/p}(\R^n).$$
Now, take $0<\mu\leq 1-n/p$. Since $W^{1,p}(\R^n)\xhookrightarrow{}C^{0,\mu}(\R^n)=\mathcal{L}^{p,p\mu+n}(\R^n)$, the same argument yields that $$H^{s,p}(\R^n)\xhookrightarrow{}\mathcal{L}^{p,\lambda_s(\mu)}(\R^n),$$ where $$\lambda_s(\mu)=(1-s)0+s(p\mu+n).$$
To have the equivalence with a H\"older space we need that $\lambda_s(\mu)>n$, so we need to restrict ourselves to $\mu\in \left(\frac{n(1-s)}{sp},1-n/p\right)$. In that case $$\mathcal{L}^{p,s(p\mu+n)}(\R^n)=C^{0,\alpha_s(\mu)}(\R^n),\,\alpha_s(\mu)=\frac{sp\mu+sn-n}{p},$$ with equivalence of norms.  The function $\mu\mapsto \alpha_s(\mu)$ is a continuous function such that $$\lim_{\mu\to \left(\frac{n(1-s)}{sp}\right)^+}\alpha_s(\mu)=0,$$ and $$\lim_{\mu\to (1-n/p)}\alpha_s(\mu)=\frac{sp(1-n/p)+sn-n}{p}=\frac{sp-n}{p}=s-\frac{n}{p}=\mu_s^*,$$ and hence it maps continuously $\left(\frac{n(1-s)}{sp},1-n/p\right)$ into $(0,\mu_s)$, so $$H^{s,p}(\R^n)\xhookrightarrow{}C^{0,\alpha}(\R^n),\,0<\alpha\leq \mu_s^*.$$

\end{enumerate}
\qed

With these results in mind and using the lifting property for the Bessel potential, we can prove the following general embedding theorem between Bessel spaces.
\begin{Teor}[Fractional Sobolev embedding theorem]\label{frset}
    Let $s,t\in (0,1)$ with $t<s$, and $1<p<\infty$ such that $(s-t)p<n$. Then, we have the embedding $$H^{s,p}(\R^n)\xhookrightarrow{}H^{t,q}(\R^n),$$
    with 
    $$q:=\frac{np}{n-(s-t)p}.$$ 
\end{Teor}
\noindent\textbf{Proof:} We have that $s>t$, so there exists $\varepsilon>0$ such that $s=t+\varepsilon$. Since $q$ is well defined by hypothesis, by {Theorem \ref{FSET}} we have that $$H^{\varepsilon,p}(\R^n)\xhookrightarrow{}L^{p_{\varepsilon}^*},$$ but note that 
$$p_{\varepsilon}^*=\frac{np}{n-\varepsilon p}=\frac{np}{n-(s-t)p}=q,$$ 
hence 
$$H^{\varepsilon,p}(\R^n)\xhookrightarrow{}L^q(\R^n),$$ 
where $L^q(\R^n)$ is identified with $H^{0,q}(\R^n)$. Then, using the fact that $\Lambda_{t}$ is an isomorphism from $H^{\varepsilon,p}(\R^n)$ to $H^{t+\varepsilon,p}(\R^n)=H^{s,p}(\R^n)$, and from $H^{0,q}(\R^n)$ to $H^{t,q}(\R^n)$, we have that 
\[\pushQED{\qed} 
H^{s,p}(\R^n)\xhookrightarrow{}H^{t,q}(\R^n)\qedhere.
\popQED\]

We have proved the 
continuous embeddings of the Bessel potential spaces
by means of the complex interpolation method, doing the
analogous as Shieh and Spector did using the Riesz fractional
gradient in \cite{ShiehSpector2015}. The fractional Sobolev embedding
states that for $1\leq p<\infty$ and $s\in (0,1)$ such
that $sp<n$, we have the continuous embedding 
$$H^{s,p}(\R^n)\xhookrightarrow{}L^{p_s^*}(\R^n),$$
where $p_s^*=np/(n-sp)$. It is well known that in the classical case,
i.e., $H^{1,p}(\R^n)=W^{1,p}(\R^n)$, replacing the Lebesgue spaces by Lorentz spaces,
we can obtain the continuous embedding (see the discussion on the introduction in \cite{cassani2013} and the references therein) $$W^{1,p}(\R^n)\xhookrightarrow{}L^{p_1^*,p},$$
and in particular since Lorentz spaces are nested with respect to the second index (as a consequence of its nature as real interpolation spaces),
$$W^{1,p}(\R^n)\xhookrightarrow{}L^{p_1^*,p}(\R^n)\xhookrightarrow{}L^{p_1^*,p_1^*}(\R^n)=L^{p_1^*}(\R^n),$$
and hence that embedding is better than the classical one. In particular, the embedding into Lorentz spaces is optimal in the 
context of rearrangement invariant spaces. 
In \cite{ShiehSpector2015} it is commented that those embeddings could be proved
in the fractional setting by means of the Riesz fractional gradient. Note that those embeddings
can not be obtained by means of interpolation since Bessel potential spaces are of complex interpolation,
the Lorentz spaces are of real interpolation, and the complex-real interpolation reiteration theorem gives us reiterated real interpolation,
and hence we would get Besov spaces. In view of this, we can not apply the techniques used before. Instead, we will use
the fractional fundamental theorem of calculus and the estimates of the Riesz potential on the Lorentz scale. For functions $u\in H^{s,p}(\R^n)$, we have that $$u=I_s*(\mathcal{R}\cdot D^su),$$ where $I_s:=\gamma_{n,s}|x|^{s-n}$ is the Riesz potential of order $0<s<n$, where $\gamma_{n,s}$ is a normalization constant (see \cite[Ch. V.1]{Stein1970}) and $\mathcal{R}$ is the Riesz transform, defined as a convolution with the Kernel $c_nx_j|x|^{-(n+1)}$ with $c_n$ a normalization constant (see \cite[Ch. III.1.2]{Stein1970}).  

The study of the mapping properties of the Riesz potential on the scale
of Lebesgue spaces was initiated by Sobolev. In 1963, O'Neil \cite{oneil1963} improved
the Hardy-Littlewood-Sobolev inequality on the Lorentz scale, proving that if $K_\alpha$ is some 
convolution operator (with suitable decay properties) of fractional order $0<\alpha<n/p$, then $K_\alpha f\in L^{p_s^*,p}(\R^n)$, for any $f\in L^{p}$.
By means of real interpolation, Peetre showed in  1966 \cite[Theorem 4.1]{peetre1966}
that the Riesz potential maps $L^p(\R^n)$ into $L^{p_s^*,p}(\R^n)$.
\begin{Teor} Let $1<p<\infty$ and $s\in (0,1)$ such that $sp<n$. Then, there exists a positive constant $C$
only depending on $s$ and $n$ such that $$\norm{I_sf}_{p_s^*,p}\leq C\norm{f}_p,$$ for any $f\in L^p(\R^n)$.
\end{Teor}

From here and the fractional fundamental theorem of calculus, the following result follows:
\begin{Teor}\label{optimosubcritico}
Let $1\leq p<\infty$ and $s\in (0,1)$ such that $sp<n$. Then, $$H^{s,p}(\R^n)\xhookrightarrow{}L^{p_s^*,q}(\R^n),\,p\leq q\leq \infty.$$
\end{Teor}
\noindent\textbf{Proof:}
Let $u\in C_c^\infty(\R^n)$. For any $v\in L^p(\R^n)$, $$\norm{I_sv}_{p_s^*,p}\leq C\norm{v}_p.$$
Now, choosing $v=\mathcal{R}\cdot \nabla^su$, since by the FFTC $I_sv=u$, $$\norm{u}_{p_s^*,p}=\norm{I_sv}_{p_s^*,p}\leq C\norm{v}_{p}\leq C'\norm{\nabla^su}_p,$$ by the $L^p$-boundedness of the Riesz transform.
Extending the result by density we get that $$H^{s,p}(\R^n)\xhookrightarrow{}L^{p_s^*,p}(\R^n),$$ and hence the result follows from the nesting on the second index of the Lorentz scale.\qed

In the critical case $sp=n$, we have established that $$H^{s,p}(\R^n)\xhookrightarrow{}L^q(\R^n),\,p\leq q<\infty,$$ and the optimal embedding $$H^{s,n/s}(\R^n)\xhookrightarrow{}\textbf{BMO}(\R^n).$$ This embedding could also be obtained by means of the fractional fundamental theorem of calculus in the spirit of Theorem \ref{optimosubcritico}. The key is the following estimate for the Riesz potential in the critical case \cite[Theorem 6.17]{duandikoetxea2001}:
\begin{Lema}\label{BMOestimate}
    Let $s\in (0,1)$ and $u\in L^{n/s}(\R^n)$. Then, there exists a constant $C$ depending only in $n$ and $s$ such that $$\norm{I_su}_{\normalfont{\textbf{BMO}}}(\R^n)\leq C\norm{u}_{n/s}.$$
\end{Lema}
From this, the embedding on \textbf{BMO} follows directly from the FFTOC.
\begin{Teor}
    Let $s\in (0,1)$. Then, the following embedding holds: $$H^{s,n/s}(\R^n)\xhookrightarrow{}\textbf{BMO}.$$
\end{Teor}
\noindent{}\textbf{Proof:}
Let $u\in C_c^\infty(\R^n)$. For any $v\in L^{n/s}(\R^n)$, $$\norm{I_sv}_{\textbf{BM0}}\leq C\norm{v}_{n/s},$$ by Lemma \ref{BMOestimate}.
Now, choosing $v=\mathcal{R}\cdot \nabla^su$, since by the FFTC $I_sv=u$, $$\norm{u}_{\textbf{BM0}}=\norm{I_sv}_{\textbf{BMO}}\leq C\norm{v}_{n/s}\leq C'\norm{\nabla^su}_{n/s},$$ by the $L^p$-boundedness of the Riesz transform for every $1<p<\infty$.
Extending the result by density we get that $$H^{s,p}(\R^n)\xhookrightarrow{}\textbf{BMO}(\R^n),$$ as we wanted to prove.\qed

Exploiting the well-known semigroup property of the Riesz potential, i.e., $I_{\alpha+\beta}=I_\alpha*I_\beta, \alpha,\beta>0$, we can give an alternative proof for Theorem \ref{frset} involving only the representation formula in \cite[Theorem 1.2]{ShiehSpector2015}, which establishes that for $u\in C_c^\infty(\R^n)$, $$D^su=I_{1-s}(Du)=D(I_ {1-s}u).$$
\begin{Teor}
    Let $0<s<t<1$ and $1<p<\infty$. Then, there exists a positive constant $C$ such that $$\norm{D^su}_{p_{t-s}^*}\leq C\norm{D^tu}_p.$$ Hence, $$H^{t,p}(\R^n)\xhookrightarrow{}H^{s,q}(\R^n),\,p\leq q\leq \frac{np}{n-(t-s)p}=p_{t-s}^*.$$
\end{Teor}
\noindent\textbf{Proof:} Since $t-s>0$, the Riesz potential $I_{t-s}$ is well defined, and by the Hardy-Littlewood-Sobolev lemma, it maps $L^p(\R^n)\to L^{p_{t-s}^*}(\R^n)$. Hence, for $u\in C_c^\infty(\R^n)$, \begin{align*}
    \norm{D^su}_{p_{t-s}^*}=\norm{I_{1-s}Du}_{p_{t-s}^*}=\norm{I_{t-s}I_{1-t}Du}_{p_{t-s}^*}\leq \norm{I_{1-t}Du}_{p}=\norm{D^tu}_{p}.
\end{align*} This implies the embedding $$H^{t,p}(\R^n)\xhookrightarrow{}H^{s,p_{t-s}^*}(\R^n).$$ The result for $p<q<p_{t-s}^*$ follows from interpolating the previous embedding with the trivial one. \qed

\subsection{Relationship Between Bessel and Gagliardo Spaces}

In the initial discussion of this section, we explored various spaces that aim to generalize Sobolev spaces to the fractional case, specifically the spaces $W^{s,p}$ (see {Examples \ref{ejemplosrealinterpol} ii)}), commonly referred to as Fractional Sobolev spaces. These spaces are crucial in the study of partial differential equations. Similar to Bessel spaces, these spaces have been introduced under different names in the literature, honoring mathematicians such as Aronszajn \cite{Aronszajn1955}, Gagliardo \cite{Gagliardo1958}, and Slobodeckij \cite{Slobodeckij1958}, who independently developed them in the late 1950s. We will refer to them as Gagliardo spaces, as they are characterized by the Gagliardo seminorm. This terminology helps distinguish them from Bessel potential spaces. The concept behind Gagliardo spaces is to extend the Hölder condition to $L^p$-functions using the Gagliardo seminorm, which measures the fractional differentiability of a function $u$ for the parameter $s \in (0,1)$. The identification of $W^{s,p}(\R^n)$ with real interpolation spaces shows a close relationship between Gagliardo and Bessel spaces, despite their apparent differences. From an interpolation perspective, they are real and complex fractional Sobolev spaces, respectively. The first significant result linking these spaces is the well-known fact that they coincide in the Hilbertian case $p=2$, a result that can be demonstrated through Harmonic Analysis and, in turn, is also a direct consequence of {Theorem \ref{PeetreTypeTheor}}.
\begin{Teor}[Coincidence of Bessel and Gagliardo spaces on the Hilbertian case]\label{Hilbertcase}
    For every $s\in (0,1)$ we have that $$H^{s,2}(\R^n)=W^{s,2}(\R^n),$$ with equivalence of norms.
\end{Teor}
\noindent\textbf{Proof:} Since $$K_{\theta,2}=\Cc_\theta,\,\theta\in (0,1),$$ when restricted to the category of Hilbert spaces. Since $L^2(\R^n)$ and $W^{1,2}(\R^n)$ are Hilbert spaces, we have that \[\pushQED{\qed} 
H^{s,2}(\R^n)=[L^2(\R^n),W^{1,2}(\R^n)]_s=\left(L^2(\R^n),W^{1,2}(\R^n)\right)_{s,2}=W^{s,2}(\R^n).\qedhere
\popQED\]
This is the only case that both spaces coincide for every fractional parameter $s$. However, they are very close to each other in the following sense.
\begin{Teor}[Contiguity of Bessel and Gagliardo spaces]\label{contiguity}
   Let $0<s_0<s<s_1<1$ and $1<p<\infty$. Then we have $$H^{s_1,p}(\R^n)\xhookrightarrow{}W^{s,p}(\R^n)\xhookrightarrow{}H^{s_0,p}(\R^n).$$ 
\end{Teor}
\noindent\textbf{Proof:} Given $\theta\in (0,1)$, by {Theorem \ref{RealReiter}} we have that $$\left(H^{s_0,p}(\R^n),H^{s_1,p}(\R^n)\right)_{\theta,p}=\left(L^p(\R^n),W^{1,p}(\R^n)\right)_{s(\theta),p}=W^{s(\theta),p}(\R^n),$$ where $s(\theta)=(1-\theta)s_0+\theta s_1$. If we choose $$\theta=\frac{s-s_0}{s_1-s_0},$$ then $s(\theta)=s$ and hence $\left(H^{s_0,p}(\R^n),H^{s_1,p}(\R^n)\right)_{\theta,p}=W^{s,p}(\R^n)$. Since $s_0<s_1$, we have that $$H^{s_0,p}(\R^n)\cap H^{s_1,p}(\R^n)=H^{s_1,p}(\R^n),\, H^{s_0,p}(\R^n)+H^{s_1,p}(\R^n)=H^{s_0,p}(\R^n),$$ so \begin{align*}H^{s_1,p}(\R^n)&=H^{s_0,p}(\R^n)\cap H^{s_1,p}(\R^n)\xhookrightarrow{}\left(H^{s_0,p}(\R^n),H^{s_1,p}(\R^n)\right)_{\theta,p}\\
&=W^{s,p}(\R^n)\xhookrightarrow{}H^{s_0,p}(\R^n)+H^{s_1,p}(\R^n)=H^{s_0,p}(\R^n).\end{align*}\qed\\

The following corollary is just a rewriting of Theorem \ref{contiguity}, and it is the usual way the result is stated in literature. 
\begin{Coro}
    Let $s\in (0,1)$ and $1<p<\infty$. For every $\varepsilon>0$ we have that $$H^{s+\varepsilon,p}(\R^n)\xhookrightarrow{}W^{s,p}(\R^n)\xhookrightarrow{}H^{s-\varepsilon,p}(\R^n).$$
\end{Coro}
Deeper into the relationship between Bessel and Gagliardo spaces, we have the following result.
\begin{Teor}\label{nesting}
    Let $1<p<\infty$ and $s\in (0,1)$. Then,$$ \begin{cases}
        W^{s,p}(\R^n)\xhookrightarrow{}H^{s,p}(\R^n),&1<p\leq 2,\\
        H^{s,p}(\R^n)\xhookrightarrow{}W^{s,p}(\R^n),&2\leq p<\infty,
    \end{cases}$$
with strict inclusions unless $p=2$.
\end{Teor}
\noindent\textbf{Proof:} Let $1<p\leq 2$. Since $L^p(\R^n)$ is of $p$-type and $W^{1,p}(\R^n)$ is isomorphic to $L^p(\R^n)$, then it is also of $p$-type, and hence by {Theorem \ref{PeetreTypeTheor}} $$W^{s,p}(\R^n)=\left(L^p(\R^n),W^{1,p}(\R^n)\right)_{s,p}\xhookrightarrow{} [L^p(\R^n),W^{1,p}(\R^n)]_s=H^{s,p}(\R^n).$$ Analogously, given $2\leq q<\infty$, since $L^q(\R^n)$ and $W^{1,q}(\R^n)$ are of $p$-type, where $1/p+1/q=1$, we have that 
\[\pushQED{\qed} 
H^{s,q}(\R^n)=[L^q(\R^n),W^{1,q}(\R^n)]_s\xhookrightarrow{}\left(L^q(\R^n),W^{1,q}(\R^n)\right)_{\theta,q}=W^{s,q}(\R^n).\qedhere
\popQED\]
This connection between Bessel and Gagliardo spaces will allow us to obtain many important results for Bessel spaces that are easier to prove in the case of Gagliardo spaces, such as compactness results, which are easily proven for Gagliardo spaces since one-sided compactness of admissible operators is preserved by real interpolation as it was finally proved by Cwikel in \cite{Cwikel1992} and by Cobos et al. in \cite{CobosKuhnSchonbek1992}. Meanwhile, the analogous case for the complex method remains an open problem; see \cite{CwikelRochberg2014} for a survey on the question.

An important problem that we have not studied is the pointwise multiplication of functions in generalized Sobolev spaces. Let $s_0, s_1 \in (0,1)$ and $1 < p_0, p_1 < \infty$, and let $X$ represent either $W$ or $H$. Given $u_j \in X^{s_j,p_j}(\R^n)$ for $j=0,1$, what can be inferred about the product $u_0 u_1$? It is anticipated that $u_0 u_1 \in X^{s,p}(\R^n)$ for appropriately chosen parameters $s$ and $p$. This classical problem has numerous applications in the modern theory of partial differential equations (PDEs) and is traditionally tackled using Littlewood-Paley theory, as well as Besov and Triebel-Lizorkin spaces. However, we have not explored this problem from the perspective of interpolation theory, as it is comprehensively addressed in the remarkable work of A. Behzadan and M. Holst \cite{BehzadanHolst2021}.

\subsection{An open problem: are $W^{1,1}(\R^n)$ and $H^{s,1}(\R^n)$ Interpolation Spaces?}\label{AppendixA}
An interesting open problem is proving or disproving whether the space $W^{1,1}(\R^n)$ is an interpolation space from a suitable couple \cite[Open problem~1.5]{ShiehSpector2018}. This comes from the question of characterizing the space of functions such that $D^su\in L^1(\R^n,\R^n)$, where $D^su$ is the Riesz fractional gradient of order $s$ \cite{ShiehSpector2015}. One would be tempted to try to use the real interpolation method; the most suitable couple is $\left(L^1(\R^n),W^{2,1}(\R^n)\right)$. Real interpolation (\cite[Theorems~17.24 and 17.30]{Leoni2017}) yields 
$$\left(L^1(\R^n),W^{2,1}(\R^n)\right)_{1/2,1}=B^1_{1,1}(\R^n),$$ 
which is a Besov space, different from $W^{1,1}(\R^n)$. The other option is to use complex interpolation of Triebel-Lizorkin spaces $F^s_{p,q}(\R^n)$. In particular $$[F_{1/2,2}^{1/2}(\R^n),F_{3/2,2}^{3/2}(\R^n)]_{1/2}=F_{1,2}^1(\R^n).$$ It is known that $F_{p,2}^s(\R^n)=H^{s,p}(\R^n)$ for $s\in\R$ and $1<p<\infty$, however, there is no known relationship for the case $p=1$. 
About the relationship between $H^{s,1}$ and $F^s_{1,2}$, we have the following observation by O. Domínguez \cite{Dominguez2025}. By simplicity, we suppose that we are working in the torus $\mathbb{T}^n$, however, there are analogous results for $\R^n$. The Hardy space $\mathcal{H}^1(\mathbb{T}^n)$ satisfies that 
$$L \log L(\mathbb{T}^n) \hookrightarrow \mathcal{H}^1(\mathbb{T}^n) \hookrightarrow L^1(\mathbb{T}^n).$$ 
This follows from the characterization of $\mathcal{H}^1(\mathbb{T}^n)$ as the space of integrable functions with an integrable Hilbert transform and the fact that the Hilbert transform of a function belonging to $L\log L(\mathbb{T}^n)$ is integrable. Applying the Fourier multiplier $m(x_j)=|x_j|^s$, $s\in \R$, with $x_j$ the j-$th$ component of $x\in\R^n$, to the embeddings, we deduce that $$H^s(L\log L)\xhookrightarrow{}F^{s}_{1,2}(\mathbb{T}^n)\xhookrightarrow{}H^{s}\left(L^1(\mathbb{T}^n)\right)=H^{s,1}(\mathbb{T}^n),$$ where we have used that the multiplier acts as a lifting and $F^0_{1,2}=\mathcal{H}^1$. Hence, the difference between $F^s_{1,2}$ and $H^{s,1}$ is a logarithmic perturbation. Note that these embeddings are optimal in the following sense. Modulo lifting, they are reduced to $$L \log L(\mathbb{T}^n) \hookrightarrow \mathcal{H}^1(\mathbb{T}^n) \hookrightarrow L^1(\mathbb{T}^n),$$ and by duality, it yields $$L^\infty(\mathbb{T}^n)\xhookrightarrow{}\textbf{BMO}(\mathbb{T}^n)\xhookrightarrow{}e^{L(\mathbb{T}^n)},$$ which is optimal due to John-Nirenberg's inequality.

Moreover, the Triebel-Lizorkin spaces coincide with the so-called "Hardy-based Sobolev spaces" $h^{k,p}(\R^n)$ (see \cite{KaltonMayborodaMitrea2007}), defined for $0<p<\infty$ and $k\in\mathbb{Z}$. It was established in \cite{MayborodaMitrea2004} that $h^{k,p}(\R^n)=F_{p,2}^k(\R^n)$ for $k\in\mathbb{Z}$ and $0<p\leq 1$, and hence by interpolation we can obtain the space $h^{1,1}(\R^n)$. However, as far as we are aware, there is no known relationship between $h^{1,1}(\R^n)$ and the classical Sobolev space $W^{1,1}(\R^n)$.\\

\section*{Acknowledgements}

This work was supported by {\it Agencia Estatal de Investigación} (Spain) through grant PID2023-151823NB-I00 and {\it Junta de Comunidades de Castilla-La Mancha} (Spain) through grant SBPLY/23/180225/000023. G.G.-S. is supported by a Doctoral Fellowship by Universidad de Castilla-La Mancha \text{2024-UNIVERS-12844-404}. We thank O. Domínguez for several conversations on the subject of this paper and in particular for his insight on the case $p=1$.

\appendix

\section{Proof of {Theorem \ref{PeetreTypeTheor}}}\label{B}

We briefly present another equivalent real interpolation method that we will require for the proof, the Method of Means or \textit{Espaces de moyennes}. The historical development of real interpolation is different: Lions and Peetre first introduced the method of means in their celebrated work \cite{LionsPeetre1964}, and then it was proved that it was in fact equivalent other methods introduced afterwards by Gagliardo, Oklander and Peetre (see \cite[Notes and Comment~3.14]{BerghLofstrom1976}). In particular, it is equivalent to the $K$-method (see \cite[Theorem~3.12.1]{BerghLofstrom1976}).\\

Let $E$ be a Banach space, $U\subseteq \R$ a $\mu$-measurable set and $1\leq p\leq \infty$. We denote by $L^p(U,E,\mu)$ the space of all strongly $\mu$-measurable $E$-valued functions $u:U\to E$ such that $$\norm{u(t)}_{L^p(U,E,\mu)}=\begin{cases}
    \left(\int_0^\infty \norm{u(t)}_E^p\,d\mu\right)^{1/p},\,1\leq p<\infty,\\
    \operatorname{sup}_{t>0}\norm{u(t)}_E,\,p=\infty
\end{cases}$$ 
is finite. Now, for a compatible couple of Banach spaces $(E_0,E_1)$, given $1\leq p_0,p_1\leq \infty$ and $\theta\in (0,1)$, we denote by $\textbf{S}\left((E_0,E_1),(p_0,p_1),\theta\right)$ the space of all $x\in E_0+E_1$ such that there exists a strongly measurable function $u:\R^+\to E_0\cap E_1$ for the $dt/t$ measure, such that 
\begin{equation}\label{eq:x}
x=\int_0^\infty \frac{u(t)}{t}\,dt,
\end{equation}with convergence in $E_0+E_1$, and such that $$\norm{t^{j-\theta}u}_{L^{p_j}(\R^+,E_j,dt/t)}<\infty,\,j=0,1.$$ 
The norm of $x$ is the infimum of
$$\operatorname{max}_{j=0,1}\left(\norm{t^{j-\theta}u}_{L^{p_j}(\R^+,E_j,dt/t)}\right)$$
over all functions $u$ such that \eqref{eq:x} holds. 

\noindent\textbf{Proof of Theorem \ref{PeetreTypeTheor}:} In fact, we will show that $$\textbf{S}\left((E_0,E_1,),(p_0,p_1),\theta)\right)\xhookrightarrow{}[E_0,E_1]_\theta\xhookrightarrow{}\textbf{S}\left((E_0,E_1,),(q_0,q_1),\theta)\right),$$ and then we will use the equivalence between the means method and the $K$-method. Let $x\in (E_0,E_1)_{\theta,p(\theta)}=\textbf{S}\left((E_0,E_1),(p_0,p_1),\theta\right)$. Then, there exists a strongly measurable function $u(t):(0,\infty)\to E_0\cap E_1$ such that $$x=\int_0^\infty \frac{u(t)}{t}\,dt,$$ with convergence in $E_0+E_1$ such that $$\norm{x}_{\textbf{S}\left((E_0,E_1),(p_0,p_1),\theta)\right)}\leq \norm{t^{-\theta}u}_{L^{p_0}\left(\R^+,E_0,dt/t)\right)}+\norm{t^{1-\theta}u}_{L^{p_1}\left(\R^+,E_1,dt/t)\right)}<\infty.$$ We fix a bounded interval $I\subset (0,\infty)$ and take $v(t)=u(t)\chi_I(t)$ and $$y=\int_0^\infty \frac{v(t)}{t}\,dt=\int_I\frac{u(t)}{t}\,dt.$$ We denote $$f(z)=\int_0^\infty t^{z-\theta}\frac{v(t)}{t}\,dt,\,z\in \C,$$ the Mellin transform of $t^{-\theta}v(t)$. The function $f$ is analytic in $\C$ and $f(\theta)=y$. Making the change of variables $t=e^{-2\pi s}$, we observe that $$f(z)=\int_{-\infty}^\infty (e^{-2\pi s})^{z-\theta}v(e^{-2\pi s})\frac{2\pi e^{-2\pi s}}{e^{-2\pi s}}\,ds=2\pi\int_{-\infty}^\infty e^{-2\pi sz}e^{2\pi s\theta}v(e^{-2\pi s})\,ds.$$ 
Now, for $r\in \R$, \begin{align*}
    f(ir)&=2\pi \int_{-\infty}^\infty e^{-2\pi i sr}e^{2\pi s\theta}v(e^{-2\pi s})\,ds=2\pi \mathcal{F}\{e^{2\pi s\theta}v(e^{-2\pi s})\}(r),\\
    f(1+ir)&=2\pi \int_{-\infty}^\infty e^{-2\pi i sr}e^{-2\pi s(1-\theta)}v(e^{-2\pi s})\,ds=2\pi \mathcal{F}\{e^{-2\pi s(1-\theta)}v(e^{-2\pi s})\}(r).
\end{align*}
Since $t^{j-\theta}u(t)\in L^{p_j}\left(\R^+,E_j,dt/t\right)$, $j=0,1,$ and the change $t=e^{-2\pi s}$ transforms the measure $dt/t$ into $-2\pi ds$, we have that $e^{-2\pi s(j-\theta)}v(e^{-2\pi s})\in L^{p_j}\left(\R^+,E_j, ds\right), j=0,1,$ and since $E_j$ is of $p_j$-type, we have $f(j+ir)\in L^{q_j}\left(\R^+,E_j,dr \right)=L^{q_j}(E_j), j=0,1$. Now, using the inequality 
$$\norm{f(\theta)}_{\theta}\leq C(\theta)\left(\norm{f(ir)}_{L^{q_0}\left(E_0\right)}+\norm{f(1+ir)}_{L^{q_1}(E_j)}\right),$$ where $C(\theta)$ is a constant depending on $\theta$, we have that 
$$\norm{y}_{\theta}\leq C(\theta)'\left(\norm{t^{-\theta}v}_{L^{p_0}\left(\R^+,E_0,dt/t)\right)}+\norm{t^{1-\theta}v}_{L^{p_1}\left(\R^+,E_1,dt/t)\right)}\right),$$ 
and taking the infimum over all possible representations $v(t)$, we get that 
$$\norm{y}_{\theta}\leq C(\theta)'\norm{y}_{\textbf{S}\left((E_0,E_1),(p_0,p_1),\theta\right)},$$ 
where the constant $C(\theta)'$ does not depend on the representation $v(t)$ or the fixed interval $I$. Now, we take some increasing family of intervals $I_m$ such that $I_m\to (0,\infty)$ as $m\to \infty$, and denote $v_m(t)=u(t)\chi_{I_m}(t)$, and $$y_m=\int_{I_m}\frac{u(t)}{t}\,dt,$$ for every $m$. Since the integral 
$$x=\int_{0}^\infty \frac{u(t)}{t}\,dt\,$$ converges in $E_0+E_1$, the sequence $\{y_m\}$ converges to $x$ in $E_0+E_1$, and by the previous calculations $$\norm{y_m}_{\theta}\leq C(\theta)'\norm{y_m}_{\textbf{S}\left((E_0,E_1),(p_0,p_1),\theta\right)},$$ for every $m$, so the sequence $\{y_m\}$ is a Cauchy sequence in $[E_0,E_1]_\theta$, and by the completeness of complex interpolation spaces it must converge to some limit $y\in [E_0,E_1]_\theta$. Since $[E_0,E_1]_\theta\xhookrightarrow{}E_0+E_1$, we have that $y=x$, and hence $$\norm{x}_\theta\leq C(\theta)'\norm{x}_{\textbf{S}\left((E_0,E_1),(p_0,p_1),\theta\right)},$$
so
$$(E_0,E_1)_{\theta,p(\theta)}=\textbf{S}\left((E_0,E_1),(p_0,p_1),\theta\right)\xhookrightarrow{}[E_0,E_1]_\theta.$$
Conversely, let $x\in [E_0,E_1]_\theta$. Then, there exists $g\in \mathfrak{F}\left(\overline{E}\right)$ such that $g(\theta)=x$. Furthermore, we can replace $g(z)$ by a function $G(z)=g(z)h(z)$ where $h(z)$ is a rapidly decreasing analytic complex-valued function such that $h(\theta)=1$, for example, taking $$G(z)=g(z)e^{-(z-\theta)^2}.$$ Then, clearly we have $G(j+it)\in L^{p_j}(E_j)$, $j=0,1.$ We define 
$$u(t)=\frac{1}{2\pi i}\int_{\gamma_a} t^{\theta-z}G(z)\,dz,$$ 
where $\gamma_a$ is the line parameterized by $\{a+is: -\infty<s<\infty\}$, $0\leq a\leq 1$. The function $u$ only takes positive values $t>0$ (indeed, it is the inverse Mellin transform of $t^{\theta}G(z)$). In the extremal cases $a=0$ and $a=1$ we obtain, respectively, $$t^{j-\theta}u(t)=\frac{1}{\pi}\int_{-\infty}^\infty t^{-is}G(j+is)\,ds,\,j=0,1.$$ Proceeding analogously as in the previous case, we could connect directly the functions $t^{j-\theta}u(t)$ with the Fourier transform of $G(j+it)$, $j=0,1,$ and since $E_j$ is of $p_j$-type, we obtain $$\norm{t^{j-\theta}u(t)}_{L^{q_j}\left(\R^+,E_j,dt/t\right)}\leq C\norm{G(j+is)}_{L^{p_j}(E_j)}\leq C' \operatorname{sup}_{s\in \R}\norm{g(j+is)}_{E_j},\,j=0,1,$$ where $C'$ is a positive constant. Now, taking $x=\theta$ in the expression of $u(t)$ yields $$u(t)=\frac{1}{2\pi i}\int_{\gamma_\theta}t^{\theta-z}G(z)\,dz=\frac{1}{2\pi}\int_{-\infty}^\infty t^{-is}G(\theta+is)\,ds$$ and taking $t=e^{2\pi r}$, $r\in \R$, $$u(e^{2\pi r })=\frac{1}{2\pi}\int_{-\infty}^\infty e^{-2\pi i rs}G(\theta+is)\,ds=(2\pi)^{-1}\mathfrak{F}\{G(\theta+is)\}(r),$$ thus $$x=G(\theta)=G(\theta+i0)=\int_{\R}\mathfrak{F}\{G(\theta+is)\}(r)\,dr=\int_{\R}2\pi u(e^{2\pi r})\,dr=\int_{0}^\infty \frac{u(t)}{t}\,dt.$$ So, we have that there exists a strongly measurable function $u(t):(0,\infty)\to E_0\cap E_1$ such that $x$ admits the representation $$x=\int_{0}^\infty \frac{u(t)}{t}\,dt,$$ with convergence in $E_0+E_1$, such that $$\norm{t^{j-\theta}u(t)}_{L^{q_j}\left(\R^+,E_j,dt/t\right)}<\infty,\,j=0,1,$$ i.e., $x\in \textbf{S}\left((E_0,E_1),(q_0,q_1),\theta\right)=(E_0,E_1)_{\theta,q(\theta)}$, hence $$(E_0,E_1)_{\theta,p(\theta)}\xhookrightarrow{}[E_0,E_1]_\theta\xhookrightarrow{}(E_0,E_1)_{\theta,q(\theta)},$$ as we wanted to prove.\qed

\section*{Conflicts of interest}

The authors declare that there are no conflicts of interest regarding the publication of this paper.

\addcontentsline{toc}{section}{References}
\bibliographystyle{plain}

\end{document}